\input amstex
\loadbold
\openup.8\jot \magnification=1200
\TagsOnRight

\def\rightwaveyarrow{\sim\!\!\sim\!\!\sim\kern-6.25pt
     \hbox{\lower.1pt\hbox{$\succ$}}}
\def\downwaveyarrow{\hbox{
     $\wr$\hbox{\kern-3pt\lower4.75pt\hbox{$\wr$}}\hbox{\kern-2.5pt
     \lower9.75pt\hbox{$\wr$}}
     \hbox{\kern-8.5pt\lower 14.5pt
     \hbox{$\curlyvee$}}}}
\def\upwaveyarrow{\hbox{
     $\curlywedge$\kern-4.60pt\lower4.0pt\hbox{$\wr$}
     \kern-6.0pt\lower11.0pt\hbox{$\wr$}
     \kern-6.0pt\lower17.5pt\hbox{$\wr$}}}

\def\bak{{\kern-1.1em }}
\def\bakk{{\kern-1.6em}}
\def\bakg{{\kern-.385em}}
\def\bakf{{\kern-.8em}}

%Greek letters
\def\alp{\alpha}		
\def\bet{\beta}		
\def\gam{\gamma}		
\def\del{\delta}		\def\Del{\Delta}
\def\eps{\varepsilon}
\def\Eta{\text H}		
\def\zet{\zeta}
\def\tet{\theta}		\def\Tet{\Theta}

\def\kap{\kappa}
\def\lam{\lambda}		\def\Lam{\Lambda}
\def\sig{\sigma}		
              
\def\vphi{\varphi}
\def\ome{\omega}		

\def\calA{{\Cal A}} \def\calB{{\Cal B}}  \def\calD{{\Cal D}}  
  \def\calI{{\Cal I}}   \def\calL{{\Cal L}}
 \def\calN{{\Cal N}}    
\def\calS{{\Cal S}}     
\def\calY{{\Cal Y}} \def\calZ{{\Cal Z}}
\def\ell{{l}}

\def\bfa{{\bold a}} \def\bfb{{\bold b}} \def\bfc{{\bold c}} \def\bfd{{\bold d}}

    \def\bfw{{\bold w}} \def\bfx{{\bold x}}
\def\bfx{{\bold x}} \def\bfy{{\bold y}} \def\bfz{{\bold z}}

 \def\bfgam{{\boldsymbol \gamma}}

 \def\bftet{{\boldsymbol \theta}}

\def\bfome{{\pmb \ome}}

\def\Lamtil{{\widetilde \Lam}}

          \def\htil{{\tilde h}}

\font\teneuf=eufm10 \font\seveneuf=eufm7 \font\fiveeuf=eufm5
\font\tenmsy=msbm10 \font\sevenmsy=msbm7 \font\fivemsy=msbm5
\font\tenmsx=msam10 \font\sixmsx=msam6 \font\fivemsx=msam5
\textfont7=\teneuf \scriptfont7=\seveneuf
\scriptscriptfont7=\fiveeuf
\textfont8=\tenmsy \scriptfont8=\sevenmsy
\scriptscriptfont8=\fivemsy
\textfont9=\tenmsx \scriptfont9=\sixmsx
\scriptscriptfont9=\fivemsx
\def\gr{\fam7 \teneuf}

\def\grB{{\gr B}}	   \def\grBhat{{\widehat \grB}}
	
\def\grD{{\gr D}}	
\def\grE{{\gr E}}	\def\gre{{\gr e}}

\def\grJ{{\gr J}}

\def\grM{{\gr M}}	\def\grm{{\gr m}}
\def\grN{{\gr N}}	\def\grn{{\gr n}}

\def\grS{{\gr S}}

  \def\dbC{{\Bbb C}}    
     
 \def\dbN{{\Bbb N}}   \def\dbQ{{\Bbb Q}} \def\dbR{{\Bbb R}}
     
 \def\dbZ{{\Bbb Z}}

\def\leaderfill{\leaders\hbox to 1em
     {\hss.\hss}\hfill}
\def\Q.E.D.{\line{\null\hfill Q.E.D.}}
\def\q.e.d.{\line{\null\hfill q.e.d.}}

\def\upchi{\hbox{\raise 2pt\hbox{$\chi$}}}
\def\upvphi{\hbox{\raise 2pt\hbox{$\vphi$}}}
\def\upgam{\hbox{\raise 2pt\hbox{$\gam$}}}
\def\VVert{\Vert\kern-1pt\vert}
\def\noint{\mathop{-\kern-10.0pt\int}}
\def\tnoint{\mathop{-\kern-8.5pt\int}}

\def\downmapsto{\downarrow\raise4.75pt\hbox{\kern-7.70pt --}}
\def\bigdownmapsto{\Big\downarrow\raise9.5pt\hbox{\kern-6pt
     --}}

\def\subsetarrow{\subset\kern-6.0pt\raise2.25pt\hbox{$\to$}}
\def\supsetarrow{\raise2.25pt\hbox{$\gets$}\kern-6.0pt\supset}
\def\arrowsim{\smash{\mathop{\to}\limits^{\lower1.5pt\hbox{$\scriptstyle
     \sim$}}}}
\def\smalltype{\font \eightrm=cmr8 \font\eightbf=cmbx8 \let\rm=\eightrm \let\bf=\eightbf \baselineskip=12pt minus 1pt \rm }
\def\le{\leq}
\def\ge{\geq}

%\TagsOnRight
\documentstyle{amsppt}
\pageheight{7.7in}
\vcorrection{-0.05in}
\NoBlackBoxes
\topmatter
\title The Hasse principle for pairs of diagonal cubic forms\endtitle
\author J\"org Br\"udern and Trevor D. Wooley$^*$\endauthor
\leftheadtext{J. Br\"udern and T. D. Wooley}
\rightheadtext{pairs of diagonal cubic forms}
\abstract By means of the Hardy-Littlewood method, we apply a new mean value 
theorem for exponential sums to confirm the truth, over the rational numbers, 
of the Hasse principle for pairs of diagonal cubic forms in thirteen or more 
variables.\endabstract
\thanks $^*$Supported in part by NSF grant DMS-010440. The authors are 
grateful to the Max Planck Institut in Bonn for its generous hospitality 
during the period in which this paper was conceived.\endthanks
\subjclass 11D72, 11L07, 11E76, 11P55\endsubjclass
\keywords Diophantine equations, exponential sums, Hardy-Littlewood method
\endkeywords
\address J\"org Br\"udern: Institut f\"ur Algebra und Zahlentheorie, 
Pfaffenwaldring 57, Universit\"at Stuttgart, D-70511 Stuttgart, Germany
\endaddress
\email bruedern\@mathematik.uni-stuttgart.de\endemail
\address Trevor D. Wooley: Department of Mathematics, University of Michigan, 
2074 East Hall, 530 Church Street, Ann Arbor, MI 48109-1043, U.S.A.\endaddress
\email wooley\@umich.edu\endemail
\endtopmatter
\document

\subhead 1. Introduction\endsubhead Early work of Lewis \cite{14} and Birch 
\cite{3, 4}, now almost a half-century old, shows that pairs of quite general 
homogeneous cubic equations possess non-trivial integral solutions whenever the
 dimension of the corresponding intersection is suitably large (modern 
refinements have reduced this permissible affine dimension to 826; see 
\cite{13}). When $s$ is a natural number, let $a_j,b_j$ $(1\le j\le s)$ be 
fixed rational integers. Then the pioneering work of Davenport and Lewis 
\cite{12} employs the circle method to show that the pair of simultaneous 
{\it diagonal} cubic equations
$$a_1x_1^3+a_2x_2^3+\ldots +a_sx_s^3=b_1x_1^3+b_2x_2^3+\ldots +b_sx_s^3=0,
\tag1.1$$
possess a non-trivial solution $\bfx \in \dbZ^s\setminus \{ {\bold 0}\} $ 
provided only that $s\ge 18$. Their analytic work was simplified by Cook 
\cite{10} and enhanced by Vaughan \cite{16}; these authors showed that the 
system (1.1) necessarily possesses non-trivial integral solutions in the cases 
$s=17$ and $s=16$, respectively. Subject to a local solubility hypothesis, 
a corresponding conclusion was obtained for $s=15$ by Baker and Br\"udern 
\cite{2}, and for $s=14$ by Br\"udern \cite{5}. Our purpose in this paper 
is the proof of a similar result that realises the sharpest conclusion 
attainable by any version of the circle method as currently envisioned, even 
if one were to be equipped with the most powerful mean value estimates for 
Weyl sums conjectured to hold.

\proclaim{Theorem 1} Suppose that $s\ge 13$, and that $a_j,b_j\in \dbZ $ 
$(1\le j\le s)$. Then the pair of equations $(1.1)$ has a non-trivial solution 
in rational integers if and only if it has a non-trivial solution in the 
$7$-adic field. In particular, the Hasse principle holds for the system $(1.1)$
 provided only that $s\ge 13$.\endproclaim

When $s\ge 13$, the conclusion of Theorem 1 confirms the Hasse principle for 
the system (1.1) in a particularly strong form: any local obstruction to 
solubility must necessarily be $7$-adic. Similar conclusions follow from the 
earlier cited work of Baker and Br\"udern \cite{2} and Br\"udern \cite{5} under
 the more stringent conditions $s\ge 15$ and $s\ge 14$, respectively.\par

The conclusion of Theorem 1 is best possible in several respects. First, when 
$s=12$, there may be arbitrarily many $p$-adic obstructions to global 
solubility. For example, let $\calS$ denote any finite set of primes $p\equiv 
1\pmod{3}$, and write $q$ for the product of all the primes in $\calS$. Choose 
any number $a \in \dbZ$ that is a cubic non-residue modulo $p$ for all $p\in 
\calS$, and consider the form
$$\Psi (x_1,\ldots ,x_6)=(x_1^3-ax_2^3)+q(x_3^3-ax_4^3)+q^2(x_5^3-ax_6^3).$$
For any $p\in \calS$, the equation $\Psi (x_1,\ldots ,x_6)=0$ has no solution 
in $\dbQ_p$ other than the trivial one, and hence the same is true of the 
pair of equations
$$\Psi (x_1,\ldots ,x_6)=\Psi (x_7,\ldots ,x_{12})=0.$$
In addition, the $7$-adic condition in the statement of Theorem 1 cannot be 
removed. Davenport and Lewis \cite{12} observed that when
$$\align \Xi (x_1,\ldots , x_5)=&\ x_1^3+2x_2^3+6x_3^3-4x_4^3,\\
\Eta (x_1,\ldots ,x_5)=&\ \ \ \ \ \ \ \ \ x_2^3+2x_3^3+4x_4^3+x_5^3,
\endalign $$
then the pair of equations in 15 variables given by
$$\align \Xi (x_1,\ldots ,x_5)+7\Xi (x_6,\ldots ,x_{10})+49\Xi (x_{11},\ldots 
,x_{15})&=0\\
\Eta (x_1,\ldots ,x_5)+7\Eta (x_6,\ldots ,x_{10})+49\Eta (x_{11},\ldots 
,x_{15})&=0\endalign $$
has no non-trivial solutions in $\dbQ_7$. In view of these examples, the state 
of knowledge concerning the local solubility of systems of the type (1.1) may 
be regarded as having been satisfactorily resolved in all essentials by 
Davenport and Lewis, and by Cook, at least when $s\ge 13$. Davenport and Lewis 
\cite{12} showed first that whenever $s\ge 16$, there are non-trivial solutions
 of (1.1) in any $p$-adic field. Later, Cook \cite{11} confirmed that such 
remains true for $13\le s\le 15$ provided only that $p\ne 7$.\par

Our proof of Theorem 1 uses analytic tools, and in particular employs the 
circle method. It is a noteworthy feature of our techniques that the method, 
when it succeeds at all, provides a lower bound for the number of integral 
solutions of (1.1) in a large box that is essentially best possible. In order 
to be more precise, when $P$ is a positive number, denote by $\calN (P)$ the 
number of integral solutions $(x_1\ldots, x_s)$ of (1.1) with $|x_j|\le P$ 
$(1\le j\le s)$. Then provided that there are solutions of (1.1) in every 
$p$-adic field, the principles underlying the Hardy-Littlewood method suggest 
that an asymptotic formula for $\calN (P)$ should hold in which the main term 
is of size $P^{s-6}$. We are able to confirm the lower bound $\calN (P)\gg 
P^{s-6}$ implicit in the latter prediction whenever the intersection (1.1) is 
in general position. This observation is made precise in the following theorem.

\proclaim{Theorem 2} Let $s$ be a natural number with $s\ge 13$. Suppose that 
$a_i,b_i\in \dbZ $ $(1\le i\le s)$ satisfy the condition that for any pair 
$(c,d) \in \dbZ^2\setminus \{(0,0)\}$, at least $s-5$ of the numbers 
$ca_j+db_j$ $(1\le j\le s)$ are non-zero. Then provided that the system $(1.1)$
 has a non-trivial $7$-adic solution, one has $\calN (P)\gg P^{s-6}$.
\endproclaim

The methods employed by earlier writers, with the exception of Cook \cite{10}, 
were not of sufficient strength to provide a lower bound for $\calN (P)$ 
attaining the order of magnitude presumed to reflect the true state of affairs.

\par The expectation discussed in the preamble to the statement of Theorem 2 
explains the presumed impossibility of a successful application of the circle 
method to establish analogues of Theorems 1 and 2 with the condition $s\ge 13$ 
relaxed to the weaker constraint $s\ge 12$. For it is inherent in applications 
of the circle method to problems involving equations of degree exceeding $2$ 
that error terms arise of size exceeding the square-root of the number of 
choices for all of the underlying variables. In the context of Theorem 2, the 
latter error term will exceed a quantity of order $P^{s/2}$, while the 
anticipated main term in the asymptotic formula for $\calN (P)$ is of order 
$P^{s-6}$. It is therefore apparent that this latter term cannot be expected 
to majorize the error term when $s\le 12$.\par

The conclusion of Theorem 2 is susceptible to some improvement. The hypotheses 
can be weakened so as to require that only seven of the numbers $ca_j+db_j$ 
$(1\le j\le s)$ be non-zero for all pairs $(c,d)\in \dbZ^2\setminus \{(0,0)\}$;
 however, the extra cases would involve us in a lengthy additional discussion 
within the circle method analysis to follow, and as it stands, Theorem 2 
suffices for our immediate purpose at hand. For a refinement of Theorem 2 along
 these lines, we refer the reader to our forthcoming communication \cite{8}.

\par In the opposite direction, we note that the lower bound recorded in the 
statement of Theorem 2 is not true without some condition on the coefficients 
of the type currently imposed. In order to see this, consider the form 
$\Psi (\bfx )$ defined by
$$\Psi (x_1,x_2,x_3,x_4)=5x_1^3+9x_2^3+10x_3^3+12x_4^3.$$
Cassels and Guy \cite{9} showed that although the equation $\Psi (\bfx )=0$ 
admits non-trivial solutions in every $p$-adic field, there are no such 
solutions in rational integers. Consequently, for any choice of coefficients 
$\bfb \in (\dbZ\setminus \{{\bold 0}\})^s$, the number of solutions $\calN (P) 
$ associated with the pair of equations
$$\Psi (x_1,x_2,x_3,x_4)=b_1x_1^3+b_2x_2^3+\ldots+ b_sx_s^3=0\tag1.2$$
is equal to the number of integral solutions $(x_5,\dots ,x_s)$ of the single 
equation $b_5x_5^3+\ldots +b_sx_s^3=0$, with $|x_i|\leq P$ $(5\le i\le s)$. 
For the system (1.2), therefore, it follows from the methods underlying 
\cite{17} that $\calN (P)\asymp P^{s-7}$ whenever $s\ge 12$. In circumstances 
in which the system (1.2) possesses non-singular $p$-adic solutions in every 
$p$-adic field, the latter is of smaller order than the prediction $\calN (P)
\asymp P^{s-6}$ consistent with the conclusion of Theorem 2 that is motivated 
by a consideration of the product of local densities. Despite the abundance of 
integral solutions of the system (1.2) for $s\ge 12$, weak approximation also 
fails. In contrast, with some additional work, our proof of Theorem 2 would 
extend to establish weak approximation for the system (1.1) without any 
alteration of the conditions currently imposed. Perhaps weak approximation 
holds for the system (1.1) with the hypotheses of Theorem 2 relaxed so as to 
require only that for any $(c,d)\in \dbZ^2\setminus \{(0,0)\}$, at least five 
of the numbers $ca_j+db_j$ $(1\le j\le s)$ are non-zero. However, in order to 
prove such a conclusion, it seems necessary first to establish that weak 
approximation holds for diagonal cubic equations in five or more variables. 
Swinnerton-Dyer \cite{15} has recently obtained such a result subject to the 
as yet unproven finiteness of the Tate-Shafarevich group for elliptic curves 
over quadratic fields.\par

This paper is organised as follows. In the next section, we announce the two 
mean value estimates that embody the key innovations of this paper; these are 
recorded in Theorems 3 and 4. Next, in section 3, we introduce a new method 
for averaging Fourier coefficients over thin sequences, and we apply it to
establish Theorem 3. Though motivated by recent work of Wooley \cite{25} and
Br\"udern, Kawada and Wooley \cite{6}, this section contains the most novel
material in this paper. In section 4, we derive Theorem 4 as well as some 
other mean value estimates that all follow from Theorem 3. Then, in section 5, 
we prepare the stage for a performance of the Hardy-Littlewood method that 
ultimately establishes Theorem 2. The minor arcs require a rather delicate 
pruning argument that depends heavily on two innovations for smooth cubic Weyl 
sums from our recent paper \cite{7}. For more detailed comments on this matter,
 the reader is directed to section 6, where the pruning is executed, and in 
particular to the comments introducing section 6. The analysis of the major 
arcs is standard, and deserves only the abbreviated discussion presented in 
section 7. In the final section, we derive Theorem 1 from Theorem 2.\par

Throughout, the letter $\eps $ will denote a sufficiently small positive 
number. We use $\ll $ and $\gg $ to denote Vinogradov's well-known notation, 
implicit constants depending at most on $\eps $, unless otherwise indicated. In
 an effort to simplify our analysis, we adopt the convention that whenever 
$\eps $ appears in a statement, then we are implicitly asserting that for each 
$\eps >0$ the statement holds for sufficiently large values of the main 
parameter. Note that the ``value'' of $\eps$ may consequently change from 
statement to statement, and hence also the dependence of implicit constants on 
$\eps$. Finally, from time to time we make use of vector notation in order to 
save space. Thus, for example, we may abbreviate $(c_1,\dots ,c_t)$ to $\bfc $.

\subhead{2. A twelfth moment of cubic Weyl sums}\endsubhead
In this section we describe the new ingredients employed in our application of 
the Hardy-Littlewood method to prove Theorem 2. The success of the method 
depends to a large extent on a new mean value estimate for cubic Weyl sums that
 we now describe. When $P$ and $R$ are real numbers with $1\leq R\leq P$, 
define the set of smooth numbers $\calA (P,R)$ by
$$\calA (P,R) = \{ n\in \dbN \cap [1,P]: p|n \text{ implies } p\le R\},$$
where, here and later, the letter $p$ is reserved to denote a prime number. The
smooth Weyl sum $h(\alp )=h(\alp ;P,R)$ central to our arguments is defined by
$$h(\alp ;P,R)=\sum_{x\in \calA (P,R)}e(\alp x^3),$$
where here and hereafter we write $e(z)$ for $e^{2\pi iz}$. An upper bound for 
the sixth moment of this sum is crucial for the discourse to follow. In order 
to make our conclusions amenable to possible future progress, we formulate the 
main estimate explicitly in terms of the sixth moment of $h(\alp )$. It is 
therefore convenient to refer to an exponent $\xi $ as {\it admissible} if, for
 each positive number $\eps $, there exists a positive number $\eta =\eta 
(\eps )$ such that, whenever $1\le R\le P^\eta$, one has the estimate
$$\int_0^1|h(\alp ;P,R)|^6\, d\alp \ll P^{3+\xi +\eps }.\tag2.1$$

\proclaim{Lemma 1} The number $\xi =(\sqrt{2833}-43)/41$ is admissible.
\endproclaim

This is the main result of \cite{22}. Since $(\sqrt{2833}-43)/41=0.2494\ldots$,
 it follows that there exist admissible exponents $\xi $ with $\xi <1/4$, a 
fact of importance to us later. The first admissible exponent smaller than 
$1/4$ was obtained by Wooley \cite{21}.

\par Next, when $a,b,c,d\in \dbZ $ and $\calB $ is a finite set of integers, we
 define the integral
$$I(a,b,c,d)=\int_0^1\!\!\int_0^1 |h(a\alp )h(b\bet )|^5\Big| \sum_{z\in \calB 
}e\big( (c\alp + d\bet )z^3\big)\Big|^2\, d\alp \, d\bet .\tag2.2$$
We may now announce our central auxiliary mean value estimate, which we prove 
in section 3.

\proclaim{Theorem 3} Suppose that $a$, $b$, $c$, $d$ are non-zero integers, 
and that $\calB \subseteq [1,P]\cap \dbZ $. Then for each admissible exponent 
$\xi $, and for each positive number $\eps $, there exists a positive number 
$\eta =\eta (\eps )$ such that, whenever $1\le R\le P^\eta $, one has
$$I(a,b,c,d)\ll P^{6+\xi +\eps }.$$
\endproclaim

If one takes $\calB =\calA (P,R)$, then the conclusion of Theorem 3 yields the 
estimate
$$\int_0^1\!\!\int_0^1 |h(a\alp )^5h(b\bet )^5h(c\alp +d\bet )^2|\, d\alp \, 
d\bet \ll P^{6+\xi +\eps }.\tag2.3$$
While this bound suffices for the applications discussed in this paper, the 
more general conclusion recorded in Theorem 3 is required in our forthcoming 
article \cite{8}. We note that previous writers would apply H\"older's 
inequality and suitable changes of variable so as to bound the left hand side 
of (2.3) in terms of factorisable double integrals of the shape
$$\int_0^1\!\!\int_0^1 |h(A\alp )h(B\bet )|^6\, d\alp \, d\bet ,\tag2.4$$
with suitable fixed integers $A$ and $B$ satisfying $AB\ne 0$. The latter 
integral may be estimated via the inequality (2.1), and thereby workers 
hitherto would derive an upper bound of the shape (2.3), but with the exponent 
$6+2\xi +\eps $ in place of $6+\xi +\eps $. Underpinning these earlier 
strategies are mean values involving two linearly independent linear forms in 
$\alp $ and $\bet $, these being reducible to the shape (2.4). In contrast, our
 approach in this paper makes crucial use of the presence within the mean value
 (2.3) of three pairwise linearly independent linear forms in $\alp $ and 
$\bet $, and we save a factor of $P^\xi $ by exploiting the extra structure 
inherent in such mean values. It is worth noting that the existence of an 
upper bound for the mean value (2.4) of order $P^{6+2\xi +\eps }$ is 
essentially equivalent to the validity of the estimate (2.1), and thus the 
strategy underlying the proof of Theorem 3 is inherently superior to that 
applied by previous authors whenever the sharpest available admissible exponent
 $\xi $ is non-zero.\par  

As another corollary of Theorem 3, we derive a more symmetric twelfth moment
estimate in section 4 below.

\proclaim{Theorem 4} Suppose that $c_i,d_i$ $(1\le i\le 3)$ are integers 
satisfying the condition
$$(c_1d_2-c_2d_1)(c_1d_3-c_3d_1)(c_2d_3-c_3d_2)\ne 0.\tag2.5$$
Write $\Lam_j=c_j\alp +d_j\bet $ $(1\le j\le s)$. Then for each admissible 
exponent $\xi $, and for each positive number $\eps $, there exists a positive
 number $\eta =\eta (\eps )$ such that, whenever $1\le R\le P^\eta $, one has 
the estimates
$$\int_0^1\!\!\int_0^1|h(\Lam_1)^5h(\Lam_2)^5h(\Lam_3)^2|\,d\alp \,d\bet \ll 
P^{6+\xi +\eps }\tag2.6$$
and
$$\int_0^1\!\!\int_0^1|h(\Lam_1)h(\Lam_2)h(\Lam_3)|^4\,d\alp \,d\bet \ll 
P^{6+\xi +\eps }.\tag2.7$$
\endproclaim

Note that the integral estimated in (2.7) has a natural interpretation as the
number of solutions of a pair of diophantine equations, an advantageous feature
absent from both (2.3) and (2.6). We remark also that conclusions analogous to 
those recorded in Theorems 3 and 4 may be derived with the cubic exponential 
sums replaced by sums of higher degree. Indeed, both the conclusions and their 
proofs are essentially identical with those presented in this paper, save that 
the admissible exponent $\xi $ herein is replaced by one depending on the 
degree in question.

\subhead{3. Averaging Fourier coefficients over thin sequences}\endsubhead
Our objective in this section is the proof of Theorem 3. We assume throughout 
that the hypotheses of the statement of Theorem 3 are satisfied. Thus, in 
particular, we may suppose that $\xi$ is admissible, and that $\eta =\eta (\eps
 )$ is a positive number sufficiently small that the estimate (2.1) holds. 
When $n\in \dbZ $, we let $r(n)$ denote the number of representations of $n$ in
 the form $n= x^3-y^3$, with $x,y\in \calB $. It follows that
$$\Big|\sum_{z\in \calB }e(\gam z^3)\Big|^2=\sum_{|n|\le P^3}r(n)e(-\gam n).
\tag3.1$$
We apply this formula to achieve a simple preliminary transformation of the 
integral $I(a,b,c,d)$ defined in (2.2). In this context, when $\ell \in \dbZ $ 
we write
$$\psi_\ell (m)=\int_0^1|h(\ell \alp )|^5e(-\alp m)d\alp .\tag3.2$$
Given $\calB \subseteq [1,P]\cap \dbZ $, the application of (3.1) within (2.2) 
leads to the relation
$$\align I(a,b,c,d)&=\sum_{|n|\le P^3} r(n)\int_0^1\!\!\int_0^1|h(a\alp )|^5
|h(b\bet )|^5e(-cn\alp )e(-dn\bet )\, d\alp \, d\bet \\
&=\sum_{|n|\le P^3}r(n)\psi_a(cn)\psi_b(dn).\endalign $$
Observe from (3.2) that $\psi_\ell (m)$ is real for any pair of integers 
$\ell $ and $m$. Then by Cauchy's inequality, we derive the basic estimate
$$I(a,b,c,d)\ll \Biggl( \sum_{|n|\le P^3}r(n)\psi_a(cn)^2\Biggr)^{1/2}\Biggl(
\sum_{|n|\le P^3}r(n)\psi_b(dn)^2\Biggr)^{1/2}.\tag3.3$$

Further progress now depends on a new method for counting integers in thin 
sequences for which certain arithmetically defined Fourier coefficients are 
abnormally large. Recent work of Wooley \cite{25} provides a framework for 
providing good estimates for the number of integers having unusually many 
representations as the sum of a fixed number of cubes. In a different 
direction, the discussion in Br\"udern, Kawada and Wooley \cite{6} supplies a 
strategy for bounding similar exceptional sets over thin sequences. Motivated
 by such arguments, we study the Fourier coefficients $\psi_\ell (km)$ for 
fixed integers $\ell $ and $k$, and in Lemma 2 below we estimate the number of 
occurrences of large values of $|\psi_\ell (kn)|$ as $n$ varies over the set 
$\calZ =\{n\in \dbZ : r(n) >0\}$. This information is then converted, in Lemma 
3, into a mean square bound for $\psi_\ell (kn)$ averaged over $\calZ $. 
Suitably positioned to bound the sums on the right hand side of (3.3), the 
proof of Theorem 3 is swiftly completed.\par

Before advancing to establish Lemma 2, we require some notation. When $\ell $ 
and $k$ are fixed integers and $T$ is a non-negative real number, we define the
 set $\calZ (T)=\calZ_{\ell ,k}(T)$ by
$$\calZ_{\ell ,k}(T)=\{ n\in \calZ \, :\, |\psi_\ell (kn)|>T\}.$$
For the remainder of this section we assume that our basic parameter $P$ is a 
large positive number, and that $\ell $ and $k$ are fixed non-zero integers.

\proclaim{Lemma 2} Whenever $\del $ is a positive number and $T\ge P^{2+\xi 
/2+ \del }$, one has the upper bound $\text{card}(\calZ (T))\ll P^{6+\xi +\eps
 }T^{-2}$.\endproclaim

\demo{Proof} We define the coefficient $\sig_m$ for each integer $m$ by means 
of the relation $\psi_\ell (m)=\sig_m|\psi_\ell (m)|$ when $\psi_\ell (m)\neq 
0$, and otherwise by putting $\sig_m=0$. Since $\calZ \subseteq [-P^3, P^3]$, 
we can define the finite exponential sum
$$K_T(\alp )=\sum_{n\in \calZ (T)}\sig_{kn}e(-kn\alp ).$$
In view of (3.2), it follows that
$$\sum_{n\in \calZ (T)}|\psi_\ell (kn)|=\int_0^1|h(\ell \alp )|^5K_T(\alp )
d\alp .\tag3.4$$
At this point, in the interest of brevity, we write $Z_T=\text{card}(\calZ (T))
$. Then the left hand side of (3.4) must exceed $TZ_T$, whence Schwarz's 
inequality yields the bound
$$TZ_T\le \left( \int_0^1 |h(\ell \alp )|^6d\alp \right)^{1/2}\left( \int_0^1
|h(\ell \alp )^4K_T(\alp )^2|d\alp \right)^{1/2}.\tag3.5$$

\par By (2.1) and a transparent change of variable, the first integral on the 
right hand side of (3.5) is $O(P^{3+\xi +\eps })$. In order to estimate the 
second integral, one first applies Weyl's differencing lemma to $|h(\ell \alp 
)|^4$ (see Lemma 2.3 of \cite{19}), and then interprets the resulting 
expression in terms of the underlying diophantine equation. Thus, one obtains
$$\int_0^1 |h(\ell \alp )^4K_T(\alp )^2|d\alp \ll P^\eps (P^3Z_T+PZ_T^2).
\tag3.6$$
For full details of this estimation, we refer the reader to Lemma 2.1 of 
Wooley \cite{24}, where a proof is described in the special case $\ell =1$ that
 readily extends to the present situation. As an alternative, we direct the 
reader to the method of proof of Lemma 5.1 of \cite{23}. Collecting together 
(3.5) and (3.6), we conclude that
$$TZ_T\ll P^{3/2 +\xi /2+\eps }(P^3Z_T+PZ_T^2)^{1/2}=P^{3+\xi /2+\eps 
}Z_T^{1/2}+P^{2+\xi /2+\eps }Z_T.$$
The proof of the lemma is completed by recalling our assumption that 
$T>P^{2+\xi /2+\del }$, where $\del $ is a positive number that we may suppose 
to exceed $2\eps $.\enddemo

\proclaim{Lemma 3} One has $\underset{n\in \calZ }\to {\sum }\psi_\ell (kn)^2
\ll P^{6+\xi +\eps }$.\endproclaim

\demo{Proof} Our discussion is facilitated by a division of the set $\calZ $ 
into various subsets. To this end, we fix a positive number $\del $ and define
$$\calY_0=\{ n\in \calZ \, :\, |\psi_\ell (kn)|\le P^{2+\xi /2+\del }\}.
\tag3.7$$
Also, when $T\ge 1$, we put $\calY (T)=\{ n\in \calZ \, :\, T<|\psi_\ell (kn)| 
\le 2T\}$. On noting the trivial upper bound $\text{card}(\calZ )\le P^2$, it 
is apparent from (3.7) that
$$\sum_{n\in \calY_0}\psi_\ell(kn)^2\le P^2(P^{2+\xi /2+\del })^2\ll P^{6+\xi 
+2\del }.\tag3.8$$
The bound $|\psi_\ell (kn)|\le P^5$, on the other hand, valid uniformly for 
$n\in \dbZ $, follows from (3.2) via the triangle inequality. A familiar 
argument involving a dyadic dissection therefore establishes that for some 
number $T$ with $P^{2+\xi /2+\del }\le T\le P^5$, one has
$$\sum_{n\in \calZ }\psi_\ell (kn)^2\ll \sum_{n\in \calY_0}\psi_\ell (kn)^2+
(\log P)\sum_{n\in \calY (T)}\psi_\ell (kn)^2.\tag3.9$$
But $\calY (T)\subseteq \calZ (T)$, and so it follows from Lemma 2 that
$$\sum_{n\in \calY (T)}\psi_\ell (kn)^2\le (2T)^2\text{card}(\calZ (T))\ll 
P^{6+\xi +\eps }.\tag3.10$$
The conclusion of Lemma 3 is obtained by substituting (3.8) and (3.10) into 
(3.9), and then taking $\del =\eps /2$.\enddemo

\proclaim{Lemma 4} One has $\underset{|n|\le P^3}\to{\sum }r(n)\psi_\ell 
(kn)^2\ll P^{6+\xi +\eps }$.\endproclaim

\demo{Proof} We begin by noting that a simple divisor argument shows that 
whenever $m$ is a non-zero integer, then $r(m)=O(P^\eps )$. Since also $r(0)\le
 P$, we find that
$$\sum_{|n|\le P^3}r(n)\psi_\ell (kn)^2\ll P\psi_\ell (0)^2+P^\eps
\sum_{n\in \calZ }\psi_\ell (kn)^2.\tag3.11$$
On recalling (3.2), moreover, it follows from a change of variable in 
combination with Schwarz's inequality that
$$\psi_\ell (0) = \int_0^1|h(\ell \alp )|^5d\alp \le \left( \int_0^1|h(\alp )
|^4d\alp \right)^{1/2}\left( \int_0^1|h(\alp )|^6d\alp \right)^{1/2}.\tag3.12$$
The first integral on the right hand side of (3.12) may be estimated by means 
of Hua's Lemma (see \cite{19, Lemma 2.5}), and the second via (2.1). Thus we 
find that
$$\psi_\ell (0)^2\ll (P^{2+\eps})(P^{3+\xi +\eps })=P^{5+\xi +2\eps }.$$
The proof of the lemma is completed by substituting the latter bound, together 
with the estimate provided by Lemma 3, into the relation (3.11).\enddemo

In order to establish Theorem 3, we have merely to apply Lemma 4 with 
$(\ell ,k)$ equal to $(a,c)$ and $(b,d)$ respectively, and then make use of the
 inequality (3.3).

\subhead{4. Some mean value estimates}\endsubhead
At this point it is convenient to explore some consequences of Theorem 3 that 
are relevant for our later proceedings. We suppose throughout this section that
 $\xi $ is admissible, and that $\eta =\eta (\eps )$ is a positive number 
sufficiently small that the estimate (2.1) holds. We begin by deriving Theorem 
4, and here we make use of the notation introduced in the statement of this 
theorem presented in section 2.

\demo{The proof of Theorem $4$} When $k$ is an integer with $1\le k\le 3$, let 
$i$ and $j$ be the integers for which $\{i,j,k\}=\{1,2,3\}$, and write
$$J_k=\int_0^1\!\!\int_0^1 |h(\Lam_i)^5h(\Lam_j)^5h(\Lam_k)^2|\,d\alp \,d\bet 
.\tag4.1$$
Then it follows from H\"older's inequality that
$$\int_0^1\!\!\int_0^1|h(\Lam_1)h(\Lam_2)h(\Lam_3)|^4d\alp \,d\bet \le 
(J_1 J_2J_3)^{1/3}.\tag4.2$$
The conclusion of Theorem 4 is immediate from the estimate $J_k=O(P^{6+\xi 
+\eps  })$ $(1\le k\le 3)$, which we now seek to establish.\par

By way of example we estimate $J_3$. Corresponding estimates for $J_1$ and 
$J_2$ follow by symmetrical arguments. We begin by observing that the 
hypotheses of Theorem 4 ensure that any two of the linear forms $\Lam_1$, 
$\Lam_2$ and $\Lam_3$ are linearly independent, whence there are non-zero 
integers $A$, $B$ and $C$, depending at most on $\bfc $ and $\bfd $, for which 
$C\Lam_3=A\Lam_1+B\Lam_2$. Making use of the periodicity (with period $1$) of 
the integrand in (4.1), and changing variables, one therefore finds that
$$\align J_3&=C^{-2}\int_0^C\!\!\int_0^C|h(\Lam_1)^5h(\Lam_2)^5h(\Lam_3)^2|
\,d\alp \,d\bet \\
&=\int_0^1\!\!\int_0^1|h(C\Lam_1)^5h(C\Lam_2)^5h(A\Lam_1+B\Lam_2)^2|\,d\alp \,
d\bet .\tag4.3\endalign $$
Now change the variables of integration from $(\alp ,\bet )$ to 
$(\Lam_1,\Lam_2)$, and observe that the resulting range of integration becomes 
a parallelogram contained in a square with sides of integral length parallel to
 the coordinate axes of the $(\Lam_1,\Lam_2)$-plane. Plainly, moreover, the 
dimensions of this square depend at most on $\bfc $ and $\bfd$. Making use 
again of the periodicity (with period $1$) of the integrand, we thus obtain the
 estimate
$$J_3\ll \int_0^1\!\!\int_0^1|h(C\Lam_1)^5h(C\Lam_2)^5h(A\Lam_1+B\Lam_2)^2|\,
d\Lam_1\,d\Lam_2.$$
The upper bound $J_3=O(P^{6+\xi+\eps })$ now follows from the consequence (2.3)
 of Theorem 3, and on making use of the corresponding symmetrical bounds for 
$J_1$ and $J_2$, the conclusion of Theorem 4 is immediate from (4.2).\enddemo

In preparation for the next lemma, we record an elementary estimate of utility 
in the arguments to follow that involve some level of combinatorial complexity.

\proclaim{Lemma 5} Let $k$ and $N$ be natural numbers, and suppose that 
$\grB \subseteq \dbC^k$ is measurable. Let $u_i(\bfz )$ $(0\le i\le N)$ be 
complex-valued functions of $\grB $. Then whenever the functions 
$|u_0(\bfz )u_j(\bfz )^N|$ $(1\le j\le N)$ are integrable on $\grB $, one has 
the upper bound
$$\int_\grB |u_0(\bfz )u_1(\bfz )\dots u_N(\bfz )|\, d\bfz \le N\max_{1\le j\le
 N}\int_\grB |u_0(\bfz )u_j(\bfz )^N|d\bfz .$$
\endproclaim

\demo{Proof} The desired conlcusion is immediate from the inequality 
$|\zet_1\zet_2\dots \zet_N|\le |\zet_1|^N+|\zet_2|^N+\dots +|\zet_N|^N$ that is
 valid for any complex numbers $\zet_i$ $(1\le i\le N)$.\enddemo

The next lemma contains (2.3) and Theorem 4 as special cases, and yet has a 
shape sufficiently general that it may be easily applied in what follows. In 
order to describe the conclusion of this lemma, we consider integers $c_j$ and 
$d_j$ with $(c_j,d_j)\ne (0,0)$ $(1\le j\le 12)$. To each pair $(c_j,d_j)$ we 
associate the linear form $\Lam_j=c_j\alp +d_j\bet $. We describe two such 
forms $\Lam_i$ and $\Lam_j$ as {\it equivalent} when there exists a non-zero 
rational number $\lam$ with $\Lam_i=\lam \Lam_j$. This notion plainly defines 
an equivalence relation on the set $\{ \Lam_1, \ldots, \Lam_{12}\}$, and we 
refer to the number of elements in the equivalence class containing the form 
$\Lam_j$ as its {\it multiplicity}. Finally, in order to promote concision, for
 each index $\ell $ we abbreviate $|h(\Lam_\ell)|$ simply to $h_\ell $. 

\proclaim{Lemma 6} In the setting described in the preamble to this lemma, 
suppose that the multiplicities of the linear forms $\Lam_1, \ldots, \Lam_{12}$
 are at most $5$. Then
$$\int_0^1\!\!\int_0^1h_1h_2\ldots h_{12}\, d\alp \,d\bet \ll P^{6+\xi+\eps}.
$$
\endproclaim

\demo{Proof} Consider the situation in which the number of equivalence classes 
amongst $\Lam_1,\ldots, \Lam_{12}$ is $t$. By relabelling indices if necessary,
 we may suppose that representatives of these equivalence classes are 
$\Lam_1, \dots ,\Lam_t$. For each index $i$, let $r_i$ denote the number of 
linear forms amongst $\Lam_1,\dots ,\Lam_{12}$ equivalent to $\Lam_i$. Then in 
view of the hypotheses of the lemma, we may relabel indices so as to ensure 
that
$$1\le r_t\le r_{t-1}\le \dots \le r_1\le 5\quad \text{and}\quad r_1+r_2+\dots 
+r_t=12.\tag4.4$$
Next, for a given index $i$ with $1\le i\le t$, consider the linear forms 
$\Lam_{\ell_j}$ $(1\le j\le r_i)$ equivalent to $\Lam_i$. Apply Lemma 5 with 
$N=r_i$, with $h_{\ell_j}$ in place of $u_j$ $(1\le j\le N)$, and with $u_0$ 
replaced by the product of those $h_\ell $ with $\Lam_\ell $ not equivalent to 
$\Lam_i$. Then it is apparent that there is no loss of generality in supposing 
that $\Lam_{\ell_j}=\Lam_i$ $(1\le j\le r_i)$. By repeating this argument for 
successive equivalence classes, moreover, we find that
$$\int_0^1\!\!\int_0^1h_1\ldots h_{12}\, d\alp \,d\bet \ll \int_0^1\!\!\int_0^1
h_1^{r_1}\ldots h_t^{r_t}\, d\alp \,d\bet .\tag4.5$$

\par A further simplification neatly sidesteps combinatorial complications. 
Let $\nu $ be a non-negative integer, and suppose that $r_{t-1}=r_t+\nu <5$. 
Then we may apply Lemma 5 with $N=\nu +2$, with $h_{t-1}$ in place of $u_i$ 
$(1\le i\le \nu +1)$ and $h_t$ in place of $u_N$, and with $u_0$ set equal to
$$h_1^{r_1}h_2^{r_2}\dots h_{t-2}^{r_{t-2}}h_{t-1}^{r_{t-1}-\nu -1}h_t^{r_t-1}.
$$
Here, and in what follows, we interpret the vanishing of any exponent as 
indicating that the associated exponential sum is deleted from the product. In 
this way we obtain an upper bound of the shape (4.5) in which the exponents 
$r_{t-1}$ and $r_t=r_{t-1}-\nu $ are replaced by $r_{t-1}+1$ and $r_t-1$, 
respectively, or else by $r_{t-1}-\nu -1$ and $r_t+\nu +1$. By relabelling if 
necessary, we derive an upper bound of the shape (4.5), subject to the 
constraints (4.4), wherein either the parameter $r_t$ is reduced, or else the 
parameter $t$ is reduced. By repeating this process, therefore, we ultimately 
arrive at a situation in which $r_{t-1}=5$, and then the constraints (4.4) 
imply that necessarily $(r_1,r_2,\dots ,r_t)=(5,5,2)$. The conclusion of the 
lemma is now immediate from (4.5) on making use of the estimate (2.6) of 
Theorem 4.\enddemo

As is often the case with applications of the circle method, it is desirable 
to have available a sharp upper bound bought with additional generating 
functions. We begin with an auxiliary lemma analogous to Theorem 4. In this 
context we take $\grm $ to be the set of real numbers $\alp \in [0,1)$ such 
that, whenever $a\in \dbZ $ and $q\in \dbN $ satisfy $(a,q)=1$ and $|q\alp -a|
\le P^{-9/4}$, then one has $q>P^{3/4}$. We then put $\grM =[0,1)\setminus 
\grm $.

\proclaim{Lemma 7} Suppose that $c_i,d_i$ $(1\le i\le 3)$ are integers 
satisfying the condition $(2.5)$. Then, in the notation employed in the 
statement of Theorem $4$, one has
$$\int_0^1\!\!\int_0^1 h_1^6h_2^6h_3^2\,d\alp \,d\bet \ll P^8.$$
\endproclaim

\demo{Proof} We observe that the argument leading from (4.1) to (4.3) reveals 
first that there are non-zero integers $A$, $B$ and $C$ for which 
$C\Lam_3=A\Lam_1+B\Lam_2$, and then via a change of variables that
$$\int_0^1\!\!\int_0^1 h_1^6h_2^6h_3^2\,d\alp \,d\bet \ll \calI_1(A,B,C),
\tag4.6$$
where we write
$$\calI_1(A,B,C)=\int_0^1\!\!\int_0^1 |h(C\alp )^6h(C\bet )^6h(A\alp +B\bet 
)^2|\,d\alp \,d\bet .\tag4.7$$
By orthogonality, the mean value $\calI_1(A,B,C)$ is bounded above by the 
number of integral solutions of the diophantine system
$$A^{-1}\sum_{i=1}^3(x_{2i-1}^3-x_{2i}^3)=B^{-1}\sum_{i=1}^3(y_{2i-1}^3-y_{2i
}^3)=C^{-1}(z_1^3-z_2^3),$$
with $1\le x_1,y_1\le P$ and $x_j,y_j,z_\ell \in \calA (P,R)$ $(2\le j\le 6,\ 
1\le \ell \le 2)$. We now introduce the classical Weyl sum
$$f(\alp ) = \sum_{1\le x\le P}e(\alp x^3),$$
and define the mean value $\calI_2(\grB )=\calI_2(\grB ;A,B,C)$, for a 
measurable set $\grB $, by
$$\calI_2(\grB )=\iint_\grB |f(C\alp )f(C\bet )h(C\alp )^5h(C\bet )^5h(A\alp 
+B\bet )^2|\, d\alp \,d\bet .\tag4.8$$
Then on applying orthogonality in combination with the triangle inequality, we 
may conclude that
$$\calI_1\le \calI_2([0,1)^2).\tag4.9$$

\par We estimate the integral on the right hand side of (4.8) by a simple
version of the circle method. By an enhanced version of Weyl's inequality (see 
\cite{17, Lemma 1}), one readily confirms that whenever $a$ is a fixed non-zero
 integer, then
$$\sup_{\tet \in \grm }|f(a\tet )|\ll P^{3/4+\eps }.\tag4.10$$
In view of the trivial upper bound $|f(\tet )|\le P$, one deduces that when 
$(\alp ,\bet )\in [0,1)^2$ and the upper bound $|f(C\alp )f(C\bet )|\ll 
P^{7/4+\eps }$ fails to hold, then necessarily $(\alp ,\bet )\in \grM^2$. 
Consequently, it follows from (4.8) and (4.9) that
$$\calI_1\ll P^{7/4+\eps }\int_0^1\!\!\int_0^1|h(C\alp )^5h(C\bet )^5h(A\alp 
+B\bet )^2|\,d\alp \,d\beta +\calI_2(\grM^2 ).\tag4.11$$
On recalling (4.7) and applying H\"older's inequality to (4.8), one finds that
$$\calI_2(\grM^2)\le \calI_1^{5/6}\Bigl( \iint_{{\grM}^2}|f(C\alp )^6f(C\bet 
)^6h(A\alp +B\bet )^2|\,d\alp \,d\bet \Bigr)^{1/6}.$$
A standard application of the Hardy-Littlewood method (see Chapter 4 of 
\cite{19}), moreover, readily confirms that whenever $a$ is a fixed non-zero 
integer, one has
$$\int_\grM |f(a\tet )|^6\,d\tet \ll P^3.$$
Thus, on making use of the trivial bound $|h(A\alp +B\bet )|\le P$, we see 
that
$$\calI_2(\grM^2)\ll \calI_1^{5/6}(P^8)^{1/6}.$$
On substituting the latter relation into (4.11) and recalling the estimate 
(2.3), we deduce that for a suitably small positive number $\del $, one has
$$\calI_1\ll P^{8-\del }+P^{4/3}\calI_1^{5/6},$$
whence $\calI_1\ll P^8$. The conclusion of the lemma is now immediate fom 
(4.6).\enddemo

With greater effort one may establish an asymptotic formula for the mean value 
recorded in the statement of Lemma 7, thereby confirming that the upper bound 
therein is of the correct order of magnitude. Were our estimate to be weaker by
 a factor of $P^\eps $, our subsequent deliberations would be greatly 
complicated.

\subhead{5. Preparing the stage for Hardy and Littlewood}\endsubhead
We are now equipped with auxiliary mean value estimates sufficient for our 
intended task, and so we return to our main concern and count integral 
solutions of the system (1.1) via the Hardy-Littlewood method. We suppose that 
the hypotheses of the statement of Theorem 2 are satisfied, so that, in 
particular, one has $s\ge 13$. With the pairs $(a_j, b_j)\in \dbZ^2\setminus 
\{(0,0)\}$ $(1\le j\le s)$, we associate both the linear forms
$$\Lam_j=a_j\alp +b_j\bet \quad (1\le j\le s),\tag5.1$$
and the two linear forms $L_1(\bftet )$ and $L_2(\bftet )$ defined for 
$\bftet \in \dbR^s$ by
$$L_1({\bftet })=\sum_{j=1}^sa_j\tet_j\quad \text{and}\quad L_2({\bftet })=
\sum_{j=1}^sb_j\tet_j.\tag5.2$$
Recall the notions of equivalence and multiplicity of linear forms from the 
preamble to Lemma 6, and extend these conventions in the natural way so as to 
apply to the set $\{ \Lam_1,\dots ,\Lam_s\}$. By the hypotheses of the 
statement of Theorem 2, one finds that for any pair $(c,d)\in \dbZ^2\setminus 
\{ (0,0)\}$, the linear form $cL_1(\bftet )+dL_2(\bftet )$ necessarily 
posesses at least $s-5$ non-zero coefficients. By choosing an appropriate 
subset $\calS $ of $\{1,\ldots, s\}$ with $\text{card}(\calS )=13$, we may 
therefore ensure that at most five of the forms $\Lam_j$ with $j\in \calS$ 
belong to the same equivalence class. Suppose that these $13$ forms fall into 
$t$ equivalence classes, and that the multiplicities of the representatives of 
these classes are $r_1,\ldots ,r_t$. In view of our earlier observations, there
 is no loss of generality in supposing that $5\ge r_1\ge r_2\ge \ldots \ge r_t$
 and $r_1+\dots +r_t=13$, and hence, in addition, that $t\ge 3$. With the aim 
of simplifying our notation, we now relabel variables in the system (1.1), and 
likewise in (5.1) and (5.2), so that the set $\calS $ becomes $\{ 1,2,\dots ,13
\}$, and so that $\Lam_1$ becomes a linear form in the first equivalence class 
counted by $r_1$, and $\Lam_2$ becomes a form in the second equivalence class 
counted by $r_2$.\par

Next, on taking suitable integral linear combinations of the equations (1.1), 
we may suppose without loss that
$$b_1=a_2=0.\tag5.3$$
Since we may suppose that $a_1b_2\ne 0$, it is now apparent that the 
simultaneous equations
$$L_1(\bftet )=L_2({\bftet })=0\tag5.4$$
possess a solution $\bftet $ with $\tet_j\ne 0$ $(1\le j\le s)$. We next apply 
the substitution $x_j\to -x_j$ for those indices $j$ with $1\le j\le s$ for 
which $\tet_j<0$. Neither the solubility of the system (1.1), nor the 
corresponding counting function $\calN (P)$, are affected by this manoeuvre, 
and yet the transformed linear system associated with (5.4) has a solution 
$\bftet $ with $\tet_j>0$ $(1\le j\le s)$. The homogeneity of the system (5.4) 
ensures, moreover, that a solution of the latter type may be chosen with 
$\bftet \in (0,1)^s$. We now fix this solution $\bftet $, and fix also $\eps $ 
to be a sufficiently small positive number, and $\eta $ to be a positive number
 sufficiently small in the context of Theorems 3 and 4 and the associated 
auxiliary mean value estimates, and so small that one has also $\eta <\tet_j<1$
 $(1\le j\le s)$. In this way, we may suppose that the solution $\bftet $ of 
the linear system (5.4) satisfies $\bftet \in (\eta ,1)^s$.\par

We are at last prepared to describe our strategy for proving Theorem 2. We take
 $P$ to be a positive number sufficiently large in terms of $\eps $, $\eta $, 
$\bfa$, $\bfb$ and $\bftet $, and we put $R = P^\eta$. On defining the 
exponential sum
$$g(\alp )=\sum_{\eta P<x\le P}e(\alp x^3)$$
and the generating functions
$$H_0(\alp ,\bet )=\prod_{j=2}^{13}h(\Lam_j)\quad \text{and}\quad H(\alp ,\bet 
)=\prod_{j=2}^sh(\Lam_j),\tag5.5$$
it follows from orthogonality that
$$\calN (P)\ge \int_0^1\!\!\int_0^1 g(\Lam_1)H(\alp ,\bet )\,d\alp \,d\bet .
\tag5.6$$
We analyse the double integral in (5.6) by means of the Hardy-Littlewood 
method. In this context, we put
$$Q=(\log P)^{1/100},\tag5.7$$
and when $a,b\in \dbZ$ and $q\in \dbN$, we define the boxes
$$\grN (q,a,b)=\{(\alp ,\bet )\in [0,1)^2\,:\,|\alp -a/q|\le QP^{-3}\text{ and 
}|\bet -b/q|\le QP^{-3}\}.$$
Our Hardy-Littlewood dissection is then defined by taking the set $\grN $ of 
{\it major arcs} to be the union of the boxes $\grN (q,a,b)$ with $0\le a,b\le 
q\le Q$ subject to $(a,b,q)=1$, and the {\it minor arcs} $\grn $ to be the 
complementary set $[0,1)^2\setminus \grN $.\par

The contribution to the integral in (5.6) arising from the major arcs $\grN $ 
satisfies the asymptotic lower bound
$$\iint_\grN g(\Lam_1)H(\alp ,\bet )\, d\alp\,d\bet \gg P^{s-6},\tag5.8$$
a fact whose confirmation is the sole objective of section 7 below. The 
corresponding contribution of the minor arcs $\grn $ is asymptotically smaller.
 Indeed, in section 6 we show that
$$\iint_\grn g(\Lam_1)H(\alp ,\bet )\,d\alp\,d\bet \ll P^{s-6}(\log P)^{
-1/140000}.\tag5.9$$
The desired conclusion $\calN (P)\gg P^{s-6}$ is immediate from (5.8) and (5.9)
 on recalling that $[0,1)^2$ is the disjoint union of $\grN $ and $\grn $.

\subhead{6. Pruning to the root}\endsubhead
Our goal in this section is the proof of the estimate (5.9). On recalling the 
definitions (5.5) and making use of the trivial bound $|h(\gam )|\le P$, it is 
apparent that the desired estimate follows directly from the following lemma, 
the proof of which will occupy us for the remainder of this section.

\proclaim{Lemma 8} Under the hypotheses prevailing in the discourse of section 
$5$, one has
$$\iint_\grn |g(\Lam_1)H_0(\alp ,\bet )|\,d\alp \,d\bet \ll P^7(\log P)^{
-1/140000}.\tag6.1$$
\endproclaim

The proof of Lemma 8 involves an unconventional pruning exercise. One gets 
started rather easily. Recall the major and minor arcs $\grM$ and $\grm $ 
introduced in the preamble to Lemma 7, and consider the auxiliary sets
$$\gre =\{(\alp ,\bet )\in \grn \,:\, \alp \in \grm \}\quad \text{and}\quad 
\grE =\{(\alp ,\bet )\in \grn \,:\, \alp \in \grM \}.\tag6.2$$
Then on recalling that $\Lam_1=a_1\alp $, one finds via two applications of 
(4.10) that
$$\sup_{(\alp ,\bet )\in \gre }|g(\Lam_1)|=\sup_{\alp \in \grm }|g(a_1\alp )|
\ll P^{3/4+\eps}.$$
But in view of the definition (5.5), the mean value of $H_0(\alp ,\bet )$ may 
be estimated by means of Lemma 6. Thus we deduce that
$$\iint_\gre |g(\Lam_1)H_0(\alp ,\bet )|\,d\alp\,d\bet \ll P^{3/4+\eps }
\int_0^1\!\!\int_0^1|H_0(\alp ,\bet )|\,d\alp\,d\bet \ll P^{27/4+\xi +\eps }.
\tag6.3$$

\par The treatment of the complementary set $\grE $ is much harder. Although 
one already has the potentially powerful information that $\alp \in \grM $, 
there is presently no such control available on $\bet $. Furthermore, there is 
only one classical Weyl sum within the product of generating functions on which
 one may hope to exercise useful control. Nonetheless, we are able to set up 
machinery with which to prune straight down to the set of narrow arcs $\grN $ 
by using two different devices from our recent work \cite{7} on cubic smooth 
Weyl sums. We are very fortunate to be able to borrow from this work, for we 
have not been successful in constructing an argument of sufficient strength 
along more conventional lines. Appropriate modifications of the aforementioned 
devices from \cite{7} are embodied in the following two lemmata.

\proclaim{Lemma 9} Let $A$ be a fixed non-zero rational number, and let $\del $
 be a fixed positive number. Then one has
$$\sup_{\lam \in \dbR }\int_\grM |g(a_1\tet )|^{2+\del }|h(A \tet +\lam )|^2\,
d\tet \ll P^{1+\del }.$$
\endproclaim

It is noteworthy, and important in our later discussion, that the bound here 
has the expected order of magnitude, uninflated by factors of $P^\eps $. The 
next lemma shares this feature.

\proclaim{Lemma 10}(i) One has
$$\int_0^1|h(\alp )|^{77/10}\,d\alp \ll P^{47/10}.$$
(ii) When $\Lam_i$ and $\Lam_j$ are inequivalent, one has
$$\iint_\grn h_i^8h_j^8\,d\alp \,d\bet \ll P^{10}Q^{-3/100}.$$
\endproclaim

We postpone the proof of these two lemmata to the end of this section, 
initiating at once the estimation of the contribution of the set $\grE $ within
 the mean value on the left hand side of (6.1). Suppose that the number of 
equivalence classes amongst $\Lam_2,\dots ,\Lam_{13}$ is $T$. By relabelling 
variables if necessary, we may suppose that representatives of these 
equivalence classes are $\Lamtil_1,\Lamtil_2,\dots ,\Lamtil_T$. For each index 
$i$, let $R_i$ denote the number of linear forms amongst $\Lam_2,\dots ,\Lam_{
13}$ equivalent to $\Lamtil_i$. Then in view of the discussion of \S5, we may 
suppose that
$$1\le R_T\le R_{T-1}\le \dots \le R_1\le 5\quad \text{and}\quad R_1+\dots 
+R_T=12.\tag6.4$$
In addition, since $\Lam_1$ has maximum multiplicity amongst $\Lam_1,\dots 
,\Lam_{13}$, and multiplicity at most $5$, we may suppose that
\item{(a)} when none of $\Lamtil_1,\dots ,\Lamtil_T$ are equivalent to 
$\Lam_1$, then necessarily $T=12$ and $R_1=R_2=\dots =R_T=1$, and
\item{(b)} when there is a linear form $\Lamtil_j$ equivalent to $\Lam_1$, then
 necessarily $R_j\le 4$.

\par Our strategy is to simplify the mean value in question using an argument 
akin to that employed in the proof of Lemma 6. First, the argument leading to 
(4.5) above in this instance shows that there is no loss of generality in 
supposing that
$$\iint_\grE |g(\Lam_1)H_0(\alp ,\bet )|\,d\alp \,d\bet \ll \iint_\grE g_1
\htil_1^{R_1}\dots \htil_T^{R_T}\,d\alp \,d\bet ,\tag6.5$$
where here, and in what follows, for each index $\ell $ we write $\htil_\ell$ 
in place of $|h(\Lamtil_\ell )|$ and $g_\ell $ in place of $|g(\Lam_\ell )|$. 
Suppose next that we are in the situation (a) above. We apply Lemma 5 with 
$N=4$, with $\htil_{3\ell -2}\htil_{3\ell -1}\htil_{3\ell }$ in place of 
$u_\ell $ $(1\le \ell \le 4)$, and with $u_0$ replaced by $g_1$. By relabelling
 indices if necessary, we obtain an upper bound of the shape (6.5) in which the
 exponent sequence $(R_1,\dots ,R_T)$ is equal to $(4,4,4)$. Now apply Lemma 5 
again with $N=2$, with $\htil_\ell $ in place of $u_\ell $ $(\ell =1,2)$, and 
with $u_0$ replaced by $g_1\htil_1^3\htil_2^3\htil_3^4$. In this way, we 
conclude that there are indices $i,j,k$ with $1<i<j<k\le 13$ for which 
$\Lam_i$, $\Lam_j$ and $\Lam_k$ are pairwise inequivalent, and
$$\iint_\grE |g(\Lam_1)H_0(\alp ,\bet )|\,d\alp \,d\bet \ll \iint_\grE 
g_1h_i^4h_j^5h_k^3\,d\alp \,d\bet .\tag6.6$$
We note for future reference the trivial observation that $\Lam_j$ is not 
equivalent to $\Lam_1$.\par

We analyse the situation (b) by applying an argument paralleling that of the 
second paragraph of the proof of Lemma 6, in this instance supposing $\nu $ to 
be a non-negative integer for which $R_{T-1}=R_T+\nu <4$, and now incorporating
 $g_1$ into the definition of $u_0$. Thus, by relabelling indices if necessary,
 we derive a bound of the shape (6.5), subject to the constraints (6.4) and 
condition (b) above, wherein $R_{T-1}=4$. The constraints (6.4) then imply 
that necessarily $(R_1,R_2,\dots ,R_T)=(5,4,3)$ or $(4,4,4)$. The latter 
circumstance may be converted to the former by means of the argument concluding
 the previous paragraph, and it is apparent that we may ensure in this process 
that $\Lamtil_2$ remains inequivalent to $\Lam_1$. In this second situation, 
therefore, we may again conclude that the bound (6.6) holds with 
$\Lam_i,\Lam_j,\Lam_k$ pairwise inequivalent, and with $\Lam_j$ not equivalent 
to $\Lam_1$.\par

Define the mean values
$$U=\iint_\grE g_1^{21/10}h_i^2h_j^{77/10}\,d\alp \,d\bet ,\quad 
V=\iint_\grE h_i^8h_k^8\,d\alp \,d\bet ,\tag6.7$$
and, when $(\ell \,m\,n)$ is a permutation of $(i\,j\,k)$,
$$W_{\ell mn}=\iint_\grE h_\ell^2h_m^6h_n^6\,d\alp \,d\bet .$$
Then a swift application of H\"older's inequality to (6.6) leads to the bound
$$\iint_\grE |g(\Lam_1)H_0(\alp ,\bet )|\,d\alp \,d\bet \ll U^{10/21}V^{1/42}
W_{ijk}^{3/84}W_{jki}^{35/84}W_{kij}^{1/21}.\tag6.8$$ 
The bound $V=O(P^{10}Q^{-3/100})$ is immediate from Lemma 10(ii), and when 
$\Lam_\ell $, $\Lam_m$ and $\Lam_n$ are pairwise inequivalent, Lemma 7 supplies
 the estimate $W_{\ell mn}=O(P^8)$. Thus we conclude from (6.8) that
$$\iint_\grE |g(\Lam_1)H_0(\alp ,\bet )|\,d\alp \,d\bet \ll P^{89/21}Q^{
-1/1400}U^{10/21}.\tag6.9$$

\par It remains to estimate the integral $U$ defined in (6.7). We recall that 
$\Lam_1=a_1\alp $, and change variables from $\bet $ to $\gam $ via the 
linear transformation $a_j\alp +b_j\bet =b_j\gam $. Note here that since 
$\Lam_1$ and $\Lam_j$ are inequivalent, then necessarily $b_j\ne 0$. Write 
$A=a_i-b_ia_j/b_j$. Then in view of the definition (6.2) of $\grE $, we may 
make use of the periodicity of the integrand to deduce that
$$U\le \int_0^1\!\!\int_\grM |g(\Lam_1)|^{21/10}|h(A\alp +b_i\gam )|^2
|h(b_j\gam )|^{77/10}\,d\alp \,d\gam \le U_1U_2,\tag6.10$$
where we write
$$U_1=\int_0^1|h(b_j\gam )|^{77/10}d\gam \quad \text{and}\quad 
U_2=\sup_{\lam \in \dbR }\int_\grM |g(a_1\alp )|^{21/10}|h(A\alp +\lam )|^2\,
d\alp .$$
An application of Lemma 10(i) reveals, via a change of variable, that 
$U_1=O(P^{47/10})$, and the bound $U_2=O(P^{11/10})$ is immediate from 
Lemma 9. Thus we find from (6.9) and (6.10) that
$$\iint_\grE |g(\Lam_1)H_0(\alp ,\bet )|\,d\alp \,d\bet \ll P^{89/21}Q^{
-1/1400}(P^{29/5})^{10/21}\ll P^7Q^{-1/1400}.$$
The conclusion of Lemma 8 now follows directly from (6.2), (6.3) and (5.7).
\vskip.1cm

We complete this section with the proofs of Lemmata 9 and 10.

\demo{The proof of Lemma $9$} Suppose that $\lam $ and $\del $ are real numbers
 with $\del >0$. Let $A$ be a fixed non-zero rational number, so that for some 
$B\in \dbZ \setminus \{ 0\}$ and $S\in \dbN $ with $(B,S)=1$, one has $A=B/S$. 
We define the modified set of major arcs $\grM^*$ by putting $\grM^*=\{ \bet 
\in [0,1):S\bet \in \grM \}$. Then a change of variable yields the relation
$$\int_\grM |g(a_1\tet )|^{2+\del }|h(A\tet +\lam )|^2\,d\tet =S\int_{\grM^*}
|g(a_1S\bet )|^{2+\del }|h(B\bet +\lam )|^2\,d\bet .\tag6.11$$
It follows from the definition of $\grM $ in the preamble to Lemma 7 that for 
each $\bet \in \grM^*$, there exist $c\in\dbZ$ and $r\in \dbN$ with 
$0\le c\le r\le P^{3/4}$, $(c,r)=1$ and $|S\bet r-c|\le P^{-9/4}$. Thus there 
exist also $a\in \dbZ $ and $q\in \dbN $ with $0\le a\le q\le SP^{3/4}$, 
$(a,q)=1$ and $|q\bet -a|\le P^{-9/4}$. We now take $\kap (q)$ to be the 
multiplicative function defined for $q\in \dbN $ by taking, for primes $p$ and
 non-negative integers $\ell $,
$$\kap (p^{3\ell })=p^{-\ell},\quad \kap (p^{3\ell +1})=2p^{-\ell -1/2},\quad 
\kap (p^{3\ell +2})=p^{-\ell -1}.$$
Then as a consequence of Theorem 4.1 and Lemmata 4.3 and 4.4 of \cite{19}, one 
has
$$\align g(a_1 S\bet) &\ll \kap (q)P(1+P^3|\bet- a/q|)^{-1}+q^{1/2+\eps }
(1+P^3|\bet -a/q|)^{1/2}\\
&\ll \kap (q)P(1+P^3|\bet -a/q|)^{-1/2}.\endalign $$
We therefore deduce from (6.11) that
$$\align \int_\grM &|g(a_1\tet )|^{2+\del }|h(A\tet +\lam )|^2\,d\tet \\
&\ll \sum_{1\le q\le SP^{3/4}}\left( \kap (q)P\right)^{2+\del }\sum^q\Sb 
a=1\\ (a,q)=1\endSb \int^\infty_{-\infty} {{\left| h\left( B(a/q+\gam )+\lam 
\right)\right|^2}\over {(1+P^3|\gam |)^{1+\del /2}}}\,d\gam .\tag6.12 
\endalign $$

\par On making use of the familiar inequality
$$\Big| \sum^q\Sb a=1\\ (a,q)=1\endSb e(a\ell/q)\Big| \le (q,\ell),$$
we find that
$$\align \sum^q\Sb a=1\\ (a,q)=1\endSb \left|h\left(B\left(a/q+\gam \right)+
\lam \right)\right|^2&=\sum_{x,y\in \calA (P,R)}\sum^q\Sb a=1\\ (a,q)=1\endSb 
e\left( (x^3-y^3)\left(B\left(a/q+\gam\right) +\lam\right)\right) \\
&\le |B|\sum_{1\le x,y\le P}(x^3-y^3,q).\endalign $$
For each natural number $q$, write $q_0$ for the cubefree part of $q$, and 
define the integer $q_3$ via the relation $q=q_0q_3^3$. Then it follows from 
the estimate (3.3) of Br\"udern and Wooley \cite{7} that whenever $1\le q\le 
P$, one has
$$\sum_{1\le x,y \le P}(x^3-y^3,q)\ll P^2q^\eps q_3.$$
In this way, we may conclude from (6.12) that
$$\align \int_\grM |g(a_1\tet )|^{2+\del }&|h(A\tet +\lam )|^2\,d\tet \\
&\ll P^{4+\del }\sum_{1\le q\le SP^{3/4}}q^\eps \kap (q)^{2+\del }q_3
\int^\infty_{-\infty}\left( 1+P^3|\gam |\right)^{-1-\del /2}\,d\gam \\
&\ll P^{1+\del }\sum_{1\le q\le SP^{3/4}}q^\eps \kap (q)^{2+\del }q_3.\tag6.13
\endalign $$
When $\del >3\eps$, moreover, the sum $\underset{q=1}\to{\overset{\infty}
\to{\sum}} q^\eps \kap(q)^{2+\del }q_3$ converges, as one readily verifies on 
recalling the definition of $\kap (q)$. The conclusion of Lemma 9 is now 
apparent from (6.13).\enddemo

\demo{The proof of Lemma $10$} The conclusion of part (i) of Lemma 10 is a 
special case of Theorem 2 of Br\"udern and Wooley \cite{7}. The proof of part 
(ii) of the lemma requires greater effort. Observe first that from Lemmata 2.2 
and 4.4 of \cite{7}, it follows easily that when $\gam $ is a real number for 
which $|h(\gam )|\ge PQ^{-1/10}$, then there exist $a\in \dbZ$ and $q\in \dbN$ 
with $1\le q\le Q^{1/3}$, $(a,q)=1$ and $|q\alp-a|\le Q^{1/3}P^{-3}$. 
Consequently, if $\Lam_k$ and $\Lam_\ell $ are inequivalent linear forms and 
$h_kh_\ell \ge P^2Q^{-1/10}$, then for $\sig =k,\ell $ there exist integers 
$d_\sig$ and $q_\sig $ with
$$1\le q_\sig \le Q^{1/3},\quad (d_\sig ,q_\sig )=1\quad \text{and}\quad 
\left| \Lam_\sig -d_\sig /q_\sig \right| \le q_\sig^{-1}Q^{1/3}P^{-3}.$$
Write $c=a_\ell b_k-a_kb_\ell $ and $\tau_c=c/|c|$, and consider the linear 
expressions $a_\ell \Lam_k-a_k\Lam_\ell $ and $b_k\Lam_\ell -b_\ell \Lam_k$. 
Then we see that in the circumstances at hand, one has $(\alp ,\bet )\in \grN 
(q,a,b)$, where
$$q=|c|q_\ell q_k,\quad a=\tau_c(b_kd_\ell q_k-b_\ell d_kq_\ell )\quad 
\text{and}\quad b=\tau_c(a_\ell d_kq_\ell -a_kd_\ell q_k).$$
It follows inter alia that when $h_kh_\ell \ge P^2Q^{-1/10}$, one necessarily 
has $(\alp ,\bet )\in \grN $. We therefore deduce that 
$$\sup_{(\alp ,\bet )\in \grn }(h_kh_\ell )\ll P^2Q^{-1/10},$$
whence
$$\iint_\grn h_k^8h_\ell^8\,d\alp\,d\bet \ll (P^2Q^{-1/10})^{3/10}\int_0^1\!\!
\int_0^1(h_kh_\ell )^{77/10}\,d\alp \,d\bet .\tag6.14$$
On making use of the first conclusion of the lemma in combination with a change
 of variables, one finds that
$$\int_0^1\!\!\int_0^1(h_kh_\ell )^{77/10}\,d\alp \,d\bet \ll \left(\int_0^1
|h(\xi )|^{77/10}\,d\xi \right)^2\ll P^{47/5},$$
and so the conclusion of the second part of the lemma follows from (6.14).
\enddemo

\subhead 7. The major arc analysis\endsubhead
We now turn our attention to the problem of estimating the contribution to the 
integral in (5.6) that arises from the major arcs $\grN $. There are relatively
 few variables involved in this integral, and our current set-up avoids various
 artifices that earlier writers have employed. For these reasons, there is no 
suitable reference available in the literature. However, the argument that we 
apply is nonetheless largely standard, and so we shall be brief.\par

First we introduce the approximants to the generating functions $g$ and $h$ on 
the major arcs $\grN $. Let
$$S(q,r)=\sum_{\ell =1}^qe(r\ell^3/q)\quad \text{and}\quad 
S_i(q,c,d)=S(q,a_ic+b_id)\quad (1\le i\le s).$$
Also, write
$$v(\tet )=\int_0^Pe(\tet \gam^3)\,d\gam \quad \text{and}\quad 
w(\tet )=\int_{\eta P}^Pe(\tet \gam^3)\,d\gam .\tag7.1$$
Finally, we mimic the convention (5.1) by associating with the pair $(a_j,b_j)$
 the linear form $\lam_j=a_j\xi +b_j\zet $ for $1\le j\le s$, and when it is 
convenient for the task at hand, we write also $v_j(\xi ,\zet )=v(\lam_j)$. 
From Lemma 8.5 of \cite{20} (see also Lemma 5.4 of \cite{18}), it follows that 
there exists a positive number $\rho $, depending at most on $\eta $, such that
 whenever $(\alp ,\bet )\in \grN (q,a,b)\subseteq \grN $, then
$$h(\Lam_j)-\rho q^{-1}S_i(q,a,b)v_i(\alp -a/q,\bet -b/q)\ll P(\log P)^{-1/2}.
\tag7.2$$
Similarly, as a consequence of Theorem 4.1 of \cite{19}, one finds that under 
the same constraints on $(\alp ,\bet )$, one has
$$g(\Lam_1)-q^{-1}S_1(q,a,b)w(a_1(\alp -a/q))\ll \log P.\tag7.3$$
Here we have made use of the hypothesis, justified by the discussion of section
 5 and recorded in (5.3), that $b_1=0$, whence in particular $\Lam_1=a_1\alp $.
 On writing
$$V(\xi ,\zet )=w(a_1\xi )\prod_{j=2}^sv(\lam_j)\quad \text{and}\quad 
U(q,a,b)=q^{-s}\prod_{j=1}^sS_i(q,a,b),\tag7.4$$
and recalling the definition (5.5), we deduce from (7.2) and (7.3) that the 
estimate
$$g(\Lam_1)H(\alp ,\bet )-\rho^{s-1}U(q,a,b)V(\alp -a/q,\bet -b/q)\ll 
P^s(\log P)^{-1/2}\tag7.5$$
holds whenever $(\alp ,\bet )\in \grN (q,a,b)\subseteq \grN $.\par

Next we introduce truncated versions of the singular integral and singular 
series, which we define respectively by
$$\grJ (X)=\iint_{\grB (X)}V(\xi ,\zet )\,d\xi \,d\zet \quad \text{and}\quad 
\grS (X)=\sum_{1\le q\le X}A(q),\tag7.6$$
in which we have written $\grB (X)$ for the box $[-XP^{-3},XP^{-3}]^2$, and 
where
$$A(q)=\underset{(c,d,q)=1}\to{\sum_{c=1}^q\sum_{d=1}^q}U(q,c,d).\tag7.7$$
The measure of the major arcs $\grN $ is $O(Q^5P^{-6})$, so that on recalling 
(5.7) and integrating over $\grN $, we infer from (7.5) that
$$\iint_\grN g(\Lam_1)H(\alp ,\bet )\,d\alp\,d\bet -\rho^{s-1}\grS (Q)\grJ (Q)
\ll P^{s-6}(\log P)^{-1/4}.\tag7.8$$

It now remains only to analyse the singular series and the singular integral 
defined, in truncated form, in (7.6). With an application in our forthcoming
article \cite{8} in mind, we study $\grS (X)$ and $\grJ (X)$ in a slightly more
 general situation than is warranted for the application at hand, and suppose 
only that for any $(c,d)\in \dbZ^2\setminus \{(0,0)\}$, at least $s-6$ of the 
numbers $ca_j+db_j$ $(1\le j\le s)$ are non-zero. In this new more general 
context, it is possible that a given linear form $\Lam_i$ may have multiplicity
 as high as six from amongst $\Lam_1,\ldots ,\Lam_{13}$. Fortunately, the 
proofs of Lemmata 12 and 13 below would be no simpler if this additional case 
were to be excluded.\par

In preparation for our discussion of the singular series, we introduce some 
additional notation and provide a simple auxiliary estimate. When $1\le j\le 
13$ and $(c,d)\in \dbZ^2\setminus \{ 0,0\}$, we define the integer 
$u_j=u_j(c,d)$ by
$$u_j=(q,ca_j+db_j).\tag7.9$$
We suppose that $\{j_1,\dots ,j_t\} \subseteq \{ 1,\dots ,13\}$ is a maximal 
set of distinct subscripts with the property that the linear forms $\Lam_{j_k}$
 are pairwise inequivalent for $1\le k\le t$. It is convenient then to define 
the integers
$$\Del =\prod_{1\le k<\ell \le t}|a_{j_k}b_{j_\ell}-a_{j_\ell }b_{j_k}|\quad 
\text{and}\quad \Xi =\Del^2\prod_{1\le k\le t}(a_{j_k},b_{j_k}).\tag7.10$$

\proclaim{Lemma 11} When $q\in \dbN $ and $(c,d)\in \dbZ^2$ satisfy the 
condition $(q,c,d)=1$, one has $u_{j_1}u_{j_2}\dots u_{j_t}|\Del q$. Moreover, 
when $v_1,\dots ,v_t$ are integers with $v_1v_2\dots v_t|\Del q$, there are 
at most $\Xi q^2(v_1\dots v_t)^{-1}$ integral pairs $(c,d)$ with $1\le c,d\le 
q$ satisfying $(c,d,q)=1$ and $u_{j_k}=v_k$ $(1\le k\le t)$.\endproclaim

\demo{Proof} Although the desired conclusions may be extracted from the 
argument of the proof of Lemma 35 of Davenport and Lewis \cite{12}, we provide 
a brief proof here for the sake of transparency of exposition. Suppose first 
that $q\in \dbN $ and $(c,d)\in \dbZ^2$ satisfy $(q,c,d)=1$. By manipulating 
appropriate linear combinations of arguments, one sees that for $1\le k<\ell
 \le t$ one has
$$(q,ca_{j_\ell}+db_{j_\ell },ca_{j_k}+db_{j_k})|(q,a_{j_\ell}b_{j_k}-a_{
j_k}b_{j_\ell })(q,c,d).$$
By hypothesis, we may suppose that $a_{j_\ell }b_{j_k}\ne a_{j_k}b_{j_\ell }$, 
and thus we deduce from (7.9) and (7.10) that
$$\prod_{1\le k<\ell \le t}(u_{j_k},u_{j_\ell })|\Del .\tag7.11$$
The desired conclusion $u_{j_1}u_{j_2}\dots u_{j_t}|\Del q$ then follows from 
the observation that $u_{j_\ell }|q$ for $1\le \ell \le t$. Next we note that 
for $1\le \ell \le t$, the number of solutions $(c,d)$, distinct modulo $v_\ell
 $, of the congruence $a_{j_\ell }c+b_{j_\ell }d\equiv 0\pmod{v_\ell }$, is 
precisely $(a_{j_\ell },b_{j_\ell },v_\ell )v_\ell $. On recalling (7.9) and 
applying the Chinese Remainder Theorem, therefore, the number of integral pairs
 $(c,d)$ satisfying $1\le c,d\le q$, $(q,c,d)=1$ and $u_{j_\ell }=v_\ell $ 
$(1\le \ell \le t)$ is at most
$$\Biggl( \prod_{1\le \ell \le t}(a_{j_\ell },b_{j_\ell })v_\ell \Biggr) 
\Biggl( q(v_1\dots v_t)^{-1}\prod_{1\le k<\ell \le t}(v_k,v_\ell )\Biggr)^2.$$
Now $(v_k,v_\ell )=(u_{j_k},u_{j_\ell })$, so on making use of (7.11) in order 
to bound the last product in this expression, we conclude from (7.10) that an 
upper bound for the number of integral pairs $(c,d)$ in question is 
$\Xi q^2(v_1\dots v_t)^{-1}$. This confirms the final conclusion of the lemma.
\enddemo

As will shortly be confirmed, the singular series $\grS $ is equal to the 
product of the $p$-adic densities of solutions. In this context we define the 
$p$-adic density $\chi_p$ by
$$\chi_p=\lim_{h\rightarrow \infty }p^{h(2-s)}M_s(p^h),\tag7.12$$
where we write $M_s(p^h)$ for the number of solutions of the system (1.1) with 
$\bfx \in (\dbZ /p^h\dbZ)^s$.

\proclaim{Lemma 12} Suppose that the linear forms $L_1(\bftet )$ and $L_2(
\bftet )$ associated with the system $(1.1)$ satisfy the condition that for any
 pair $(c,d) \in\dbZ^2\setminus \{(0,0)\}$, the linear form $cL_1(\bftet 
)+dL_2(\bftet )$ contains at least $s-6$ non-zero coefficients. Then the limit 
$\grS = \lim\limits_{X\to \infty }\grS (X)$ exists, and one has
$$\grS -\grS (X)\ll X^{-1/4}.\tag7.13$$
Moreover, the Euler product $\prod_p \chi_p$ converges absolutely to $\grS $, 
and one has $\chi_p =0$ if and only if the system $(1.1)$ has no non-trivial 
solution in $\dbQ_p$. Finally, when the system $(1.1)$ possesses a non-trivial 
solution in $\dbQ_p$ for every prime number $p$, one has $\grS \gg 1$.
\endproclaim

\demo{Proof} We establish the upper bound
$$A(q)\ll q^{\eps-4/3}.\tag7.14$$
In view of (7.6), the estimate (7.14) not only confirms (7.13) but also shows 
that the limit $\grS =\lim\limits_{X\to \infty }\grS (X)$ exists. The proof of 
the remaining conclusions of the lemma follow by the theory familiar to 
practitioners of the circle method (see, for example, Section 2.6 of \cite{19},
 or Section 10 of Davenport and Lewis \cite{12}). From (7.7) and (7.14) one 
finds that
$$p^{h(2-s)}M_s(p^h)=\sum_{\ell =0}^{p^h}A(p^\ell )=1+O(p^{\eps -4/3}),$$
so that the definition (7.12) yields the estimate $\chi_p=1+O(p^{-5/4})$. It 
follows that the Euler product $\prod_p\chi_p$ converges absolutely. But from 
(7.7) we see that $A(q)$ is a multiplicative function of $q$, and so we see 
from (7.6) that indeed $\grS $ is equal to the aforementioned Euler product. 
We may choose $p_0$ large enough so that $\chi_p\ge 1-p^{-6/5}$ for $p>p_0$, 
and then it follows that
$$\grS \gg \prod_{p\le p_0}\chi_p.\tag7.15$$

\par In circumstances wherein the system (1.1) has no non-trivial solution in 
$\dbQ_p$, it fails to possess a non-trivial solution in $\dbZ /p^h\dbZ $ for 
any $h$, and thus it follows from (7.12) that $\chi_p=0$. On the other hand, 
when the system (1.1) possesses a non-trivial solution in $\dbQ_p$, the 
argument of the proof of the Corollary to Theorem 1 of Davenport and Lewis 
\cite{12} (see the end of section 5 of the latter paper) shows that the system 
(1.1) has a non-singular solution in $\dbQ_p$. An argument employing Hensel's 
Lemma (as in Lemma 6.7 of \cite{20}, for example) then shows that $M_s(p^h)\gg 
p^{h(2-s)}$ for large enough values of $h$, whence (7.12) shows that 
$\chi_p>0$. Thus $\chi_p=0$ if and only if the system (1.1) has no non-trivial 
solution in $\dbQ_p$, and by (7.15) one has $\grS \gg 1$ unless the system 
(1.1) fails to possess a non-trivial solution in $\dbQ_p$ for some prime $p$.

\par It remains to establish (7.14). First, by relabelling indices if 
necessary, the hypotheses of the lemma permit the assumption that the maximum 
multiplicity of any of the forms $\Lam_1,\dots ,\Lam_{13}$ is six. By Theorem 
4.2 of \cite{19}, whenever $(q,r)=1$ one has $S(q,r)\ll q^{2/3}$. Thus, on 
recalling the definition (7.9), one finds that
$$q^{-1}S(q,ca_j+db_j)\ll u_j^{1/3}q^{-1/3}.$$
Consequently, using trivial estimates for factors in the definition (7.4) of 
$U(q,c,d)$ with $j>13$, we deduce that
$$U(q,c,d)\ll q^{-13/3}(u_1 u_2 \ldots u_{13})^{1/3}.$$
We note now that when $\Lam_\ell $ and $\Lam_k$ are equivalent linear 
forms, it is a consequence of (7.9) that $u_\ell \ll u_k\ll u_\ell $. Recall 
the notation defined in the preamble to Lemma 11, and suppose that for $1\le k
\le t$, the linear form $\Lam_{j_k}$ has multiplicity $r_k$ amongst $\Lam_1,
\dots ,\Lam_{13}$. Then in view of (7.7), it follows from Lemma 11 that
$$\align A(q)&\ll q^{-13/3}\underset{(c,d,q)=1}\to{\sum_{c=1}^q\sum_{d=1}^q}
(u_1u_2\ldots u_{13})^{1/3}\\
&\ll q^{-13/3}\sum \Sb v_1,\dots ,v_t\\ v_1\ldots v_t|\Del q\endSb {{\Xi q^2}
\over {v_1\ldots v_t}}(v_1^{r_1}\ldots v_t^{r_t})^{1/3}.\endalign $$
Observe that since the linear forms $\Lam_{j_1},\dots ,\Lam_{j_t}$ are 
pairwise inequivalent, the integer $\Del $ is non-zero, and further, the 
integers $\Del $ and $\Xi $ are bounded purely in terms of the coefficients 
$\bfa $ and $\bfb $. We are permitted to assume that $r_j\le 6$ for $1\le j\le 
t$, so on using an elementary bound for the divisor function, we conclude that
$$A(q)\ll q^{-7/3}\sum \Sb v_1,\dots ,v_t\\ v_1\dots v_t|\Del q\endSb v_1\dots 
v_t\ll q^{\eps - 4/3},$$
as claimed in (7.14). This completes the proof of the lemma.\enddemo

We now turn our attention to the truncated singular integral $\grJ (X)$. The 
analysis here is very straightforward, but ironically, the simplicity of our 
approach prevents any convenient reference to the literature.

\proclaim{Lemma 13} Under the same hypotheses as in the statement of Lemma 
$12$, the limit $\grJ =\lim\limits_{X\to \infty }\grJ (X)$ exists, and one has
$$\grJ -\grJ (X)\ll P^{s-6}X^{-1}.\tag7.16$$
Moreover, one has $\grJ \gg P^{s-6}$.\endproclaim

\demo{Proof} We begin by considering two inequivalent forms $\Lam_i$ and 
$\Lam_j$. When $T$ is a positive number, write $\grB (T)$ for the box 
$[-T,T]^2$, and $\grBhat (T)$ for the complementary set $\dbR^2\setminus \grB 
(T)$. Consider now a positive number $Y$ and suppose that $(\Lam_i,\Lam_j)\in 
\grB (Y)$. The latter assertion is equivalent to the statement that
$$a_i\alp +b_i\bet \in [-Y,Y]\quad \text{and}\quad a_j\alp +b_j\bet \in 
[-Y,Y],$$
so that on taking suitable linear combinations of these forms, one obtains
$$(a_ib_j-a_jb_i)\alp \in |b_j|\, [-Y,Y]+|b_i|\, [-Y,Y]$$
and
$$(a_jb_i-a_ib_j)\bet \in |a_j|\, [-Y,Y]+|a_i|\, [-Y,Y],$$
the sums of intervals being interpreted set theoretically. The integer 
$a_ib_j-a_jb_i$ is non-zero, so that if we write
$$\Tet^{-1}=\max_{1\le i<j\le s}\{ |a_i|+|a_j|\,,\,|b_i|+|b_j|\} ,$$
then we may conclude that $(\alp ,\bet )\in \grB (\Tet^{-1}Y)$. On making a 
transparent change of variables, it follows from this discussion that when 
$\Lam_i$ and $\Lam_j$ are inequivalent, one has
$$\iint_{\grBhat (XP^{-3})}|v(\lam_i)v(\lam_j)|^6\,d\xi \,d\zet \ll \iint_{
\grBhat (\Tet XP^{-3})}|v(\lam_i)v(\lam_j)|^6\,d\lam_i\, d\lam_j.\tag7.17$$

\par We now recall the estimate $v(\gam )\ll P(1+P^3|\gam |)^{-1/3}$ that 
follows, for example, from Theorem 7.3 of \cite{19}. In the situation at hand, 
we may suppose that none of the forms $\Lam_1,\dots ,\Lam_{13}$ has 
multiplicity exceeding six. Hence, following a suitable relabelling of indices,
 we may temporarily suppose that for $1\le k\le 6$, the forms $\Lam_{2k}$, 
$\Lam_{2k+1}$ are inequivalent. Then on integrating the elementary inequality
$$|v(\lam_2)v(\lam_3)\ldots v(\lam_{13})|\le \sum_{k=1}^6|v(\lam_{2k})v(
\lam_{2k+1})|^6$$
and applying the observation (7.17), we deduce that
$$\align \iint_{\grBhat (XP^{-3})}&|v(\lam_2)v(\lam_3)\dots v(\lam_{13})|\,
d\xi \,d\zet \\
&\ll P^{12}\iint_{\grBhat (\Tet XP^{-3})}(1+P^3|\lam |)^{-2}(1+
P^3|\mu |)^{-2}\,d\lam \,d\mu \\
&\ll P^6X^{-1}.\endalign $$
Finally, on using trivial bounds for $w(a_1\xi )$ and $v(\lam_j)$ when $j>13$, 
we conclude that
$$\iint_{\grBhat (XP^{-3})}|V(\xi ,\zet )|\,d\xi \,d\zet \ll P^{s-6}X^{-1}.
\tag7.18$$
In particular, it follows from (7.6) that the singular integral $\grJ =\lim 
\limits_{X\rightarrow \infty }\grJ (X)$ exists, and (7.18) provides the 
desired estimate (7.16).\par

In order to evaluate $\grJ $ we follow the familiar routine based on the use 
of Fourier's integral theorem. From (7.1) and (7.6), we see that
$$\grJ =\int_{-\infty }^\infty \!\int_{-\infty }^\infty \!\int_\grD e(\xi 
\calL_1(\bfgam )+\zet \calL_2(\bfgam ))\, d\bfgam \, d\xi \, d\zet ,$$
where we write
$$\calL_1(\bfgam )=\sum_{i=1}^sa_i\gam_i^3\quad \text{and}\quad \calL_2(
\bfgam )=\sum_{i=1}^sb_i\gam_i^3,$$
and where $\calD $ denotes the box $[\eta P,P]\times [0,P]^{s-1}$. Put 
$\mu =P^3\xi $ and $\nu =P^3\zet $, and substitute $\ome_i=(P^{-1}\gam_i)^3$ 
for $1\le i\le s$. Then with these changes of variables we discover that
$$\grJ =3^{-s}P^{s-6}\int_{-\infty}^\infty \!\int_{-\infty}^\infty \int_{
\grD'}{{e(\mu L_1(\bfome )+\nu L_2(\bfome ))}\over {(\ome_1\dots \ome_s)^{2/3
}}}\,d\bfome \,d\mu\, d\nu ,$$
where $L_1(\bfome )$ and $L_2(\bfome )$ are defined in (5.2), and where 
$\grD'=[\eta^3,1]\times [0,1]^{s-1}$. The discussion of section 5 ensures that 
the equations $L_1(\bfome )=L_2(\bfome )=0$ define an $(s-2)$-dimensional 
linear space that passes through the point $\bftet $ lying in the interior of 
$\grD'$. Recall from (5.3) that $b_1=a_2=0$, whence from (5.2),
$$\ome_1=a_1^{-1}\Bigl( L_1(\bfome )-\sum_{j=3}^sa_j\ome_j\Bigr) \quad 
\text{and}\quad \ome_2=b_2^{-1}\Bigl( L_2(\bfome )-\sum_{j=3}^sb_j\ome_j\Bigr) 
.$$
Then on making a change of variables and applying Fourier's integral formula 
twice, in the shape
$$\lim_{\lam \rightarrow \infty }\int_{-T}^T\!\!\int_{-\lam }^\lam V(t)e(t\ome 
)\,d\ome \,dt =V(0),$$
we obtain the relation
$$\grJ \gg P^{s-6}\int_{\grD''}(\ome_1\dots \ome_s)^{-2/3}\,d\ome_3\,d\ome_4
\dots d\ome_s.\tag7.19$$
Here, we define the coordinates $\ome_1$ and $\ome_2$ by
$$\ome_1=-a_1^{-1}\sum_{j=3}^sa_j\ome_j\quad \text{and}\quad \ome_2=-b_2^{-1}
\sum_{j=3}^sb_j\ome_j,$$
and we write $\grD''$ for the set of $(s-2)$-tuples $(\ome_3,\ome_4,\dots 
,\ome_s)\in [0,1]^{s-2}$ for which the $s$-tuple $(\ome_1,\dots ,\ome_s)$ lies 
in $\grD'$. Notice that the point $(\tet_3,\dots ,\tet_s)$ necessarily lies in 
the interior of the polytope $\grD''$, whence $\grD''$ has positive volume. The
 latter observation ensures that the integral on the right hand side of (7.19) 
is positive. Since, plainly, the latter integral is independent of $P$, we may 
conclude that $\grJ \gg P^{s-6}$, and this completes the proof of the lemma.
\enddemo

The proof of Theorem 2 is now swiftly completed. By (7.8) and Lemmata 12 and 
13, we find that
$$\iint_\grN g(\Lam_1)H(\alp ,\bet )\,d\alp\,d\bet -\rho^{s-1}\grS \grJ \ll 
P^{s-6}Q^{-1/4},$$
so that in view of the estimate (5.9) we may conclude that
$$\int_0^1\!\!\int_0^1g(\Lam_1)H(\alp ,\bet )\,d\alp\,d\bet =\rho^{s-1}\grS \grJ +O(
P^{s-6}(\log P)^{-1/140000}).$$
From Lemmata 12 and 13, moreover, it is apparent that $\rho^{s-1}\grS \grJ 
\gg P^{s-6}$ provided only that the system (1.1) has non-trivial solutions in 
$\dbQ_p$ for every prime $p$. But in such circumstances, the lower bound (5.6) 
ensures that $\calN (P)\gg P^{s-6}$. In view of the discussion on $p$-adic 
solubility prior to the statement of Theorem 2, solubility over $\dbQ_p$ is 
already assured when $p\ne 7$, and the conclusion of Theorem 2 follows 
immediately.

\subhead 8. Le coup de gr\^ace\endsubhead
The theme of this concluding section is the proof of Theorem 1. Needless to 
say, if Theorem 2 is applicable to the system (1.1), then there is nothing 
further to discuss. Thus we may suppose that there exists a pair $(c,d)\in 
\dbZ^2\setminus \{(0,0)\}$ with the property that at most $s-6$ of the numbers 
$ca_j+db_j$ $(1\le j\le s)$ are non-zero. By taking suitable rational linear 
combinations of the two equations defining (1.1), it is apparent that there is 
no loss of generality in supposing that the system (1.1) takes the shape
$$a_1x_1^3+\ldots +a_sx_s^3=b_{t+1}x_{t+1}^3+\ldots +b_sx_s^3=0,\tag8.1$$
where $s\ge 13$ and $t\ge 6$. We now recall a conclusion of R. Baker 
concerning the solubility of diagonal cubic equations.

\proclaim{Lemma 14} Whenever $r\ge 7$ and $c_1,\dots ,c_r$ are rational 
integers, the equation $c_1x_1^3+\ldots +c_rx_r^3=0$ possesses a non-trivial 
integral solution.\endproclaim

\demo{Proof} On setting $x_i=0$ for $i>7$, the desired conclusion follows from 
\cite{1}.\enddemo

Let us return to the system (8.1). If one has $t\ge 7$, then it follows from 
Lemma 14 that the equation $a_1y_1^3+\ldots +a_ty_t^3=0$ possesses a 
non-trivial integral solution $\bfy =\bfz $, and thus the system (8.1) has the 
non-trivial integral solution $\bfx =(z_1,\dots ,z_t,0,\dots ,0)$. We are 
therefore left to ponder the situation in which $t=6$ and $s\ge 13$. In view 
of Lemma 14, the equation $b_7y_7^3+\ldots +b_sy_s^3=0$ possesses a non-trivial
 integral solution $(y_7,\dots ,y_s)=(z_7,\dots ,z_s)$. We put $A=a_7z_7^3+
\dots +a_sz_s^3$, and consider the equation $Ay_0^3+a_1y_1^3+\ldots 
+a_6y_6^3=0$. This equation possesses a non-trivial solution $\bfy =\bfw $, 
again by Lemma 14, and so the system (8.1) in this instance has the non-trivial
 integral solution
$$(x_1,\dots ,x_s)=(w_1,\ldots ,w_6,w_0z_7,\dots ,w_0z_s).$$
We therefore conclude that when $s\ge 13$ and Theorem 2 fails to deliver 
the Hasse principle for the system (1.1), this system nonetheless possesses 
non-trivial integral solutions. This completes the proof of Theorem 1.


\Refs
\widestnumber \no{15}

\ref \no1 \by R. C. Baker\paper Diagonal cubic equations II\jour Acta Arith. 
\vol 53\yr 1989\pages 217--250\endref

\ref \no2 \by R. C. Baker and J. Br\"udern\paper On pairs of additive cubic 
equations\jour J. Reine Angew. Math.\vol 391\yr 1988\pages 157--180\endref

\ref \no3 \by B. J. Birch\paper Homogeneous forms of odd degree in a large 
number of variables\jour Mathematika\vol 4\yr 1957\pages 102--105\endref

\ref \no4 \by B. J. Birch\paper Forms in many variables\jour Proc. Roy. Soc. 
Ser. A\vol 265\yr 1961\pages 245--263\endref

\ref \no5 \by J. Br\"udern\paper On pairs of diagonal cubic forms\jour Proc. 
London Math. Soc. (3)\vol 61\yr 1990\pages 273--343\endref

\ref \no6 \by J. Br\"udern, K. Kawada and T. D. Wooley\paper Additive 
representation in thin sequences, I: Waring's problem for cubes\jour 
Ann. Sci. \'Ecole Norm. Sup. (4)\vol 34\yr 2001\pages 471--501\endref

\ref \no7 \by J. Br\"udern and T. D. Wooley\paper On Waring's problem for 
cubes and smooth Weyl sums\jour Proc. London Math. Soc. (3)\vol 82\yr 2001
\pages 89--109\endref

\ref \no8 \by J. Br\"udern and T. D. Wooley\paper The density of integral 
solutions for pairs of diagonal cubic equations\yr preprint\endref

\ref \no9 \by J. W. S. Cassels and M. J. T. Guy\paper On the Hasse principle 
for cubic surfaces\jour Mathematika\vol 13\yr 1966\pages 111--120\endref

\ref \no10 \by R. J. Cook\paper Pairs of additive equations\jour Michigan 
Math. J. \vol 19\yr 1972\pages 325--331\endref

\ref \no11 \by R. J. Cook\paper Pairs of additive congruences: cubic 
congruences\jour Mathematika\vol 32\yr 1985\pages 286--300\endref

\ref \no12 \by H. Davenport and D. J. Lewis\paper Cubic equations of additive 
type\jour Philos. Trans. Roy. Soc. London Ser. A\vol 261\yr 1966\pages 
97--136\endref

\ref \no13 \by R. Dietmann and T. D. Wooley\paper Pairs of cubic forms in many 
variables\jour Acta Arith. \vol 110\yr 2003\pages 125--140\endref

\ref \no14 \by D. J. Lewis\paper Cubic forms over algebraic number fields
\jour Mathematika\vol 4\yr 1957\pages 97--101\endref

\ref \no15 \by H. P. F. Swinnerton-Dyer\paper The solubility of diagonal cubic 
surfaces\jour Ann. Sci. \'Ecole Norm. Sup. (4)\vol 34\yr 2001\pages 891--912
\endref

\ref \no16 \by R. C. Vaughan\paper On pairs of additive cubic equations
\jour Proc. London Math. Soc. (3)\vol 34\yr 1977\pages 354--364\endref

\ref \no17 \by R. C. Vaughan\paper On Waring's problem for cubes\jour J. Reine
 Angew. Math. \vol 365\yr 1986\pages 122--170\endref

\ref \no18 \by R. C. Vaughan\paper A new iterative method in Waring's 
problem\jour Acta Math.\vol 162\yr 1989\pages 1--71\endref

\ref \no19 \by R. C. Vaughan\book The Hardy-Littlewood Method, second edition
\publ Cambridge University Press\yr 1997\endref

\ref \no20 \by T. D. Wooley\paper On simultaneous additive equations, II\jour 
J. Reine Angew. Math.\vol 419\yr 1991\pages 141--198\endref

\ref \no21 \by T. D. Wooley\paper Breaking classical convexity in Waring's 
problem: sums of cubes and quasi-diagonal behaviour\jour Invent. Math. 
\vol 122\yr 1995\pages 421--451\endref

\ref \no22 \by T. D. Wooley\paper Sums of three cubes\jour Mathematika\vol 47
\yr 2000\pages 53--61\endref

\ref \no23 \by T. D. Wooley\paper Slim exceptional sets for sums of cubes
\jour Canad. J. Math.\vol 54 \yr 2002\pages 417--448\endref

\ref \no24 \by T. D. Wooley\paper Slim exceptional sets and the asymptotic 
formula in Waring's problem\jour Math. Proc. Cambridge Philos. Soc.\vol 134
\yr 2003\pages 193--206\endref

\ref \no25 \by T. D. Wooley\paper Slim exceptional sets for stout 
representation problems\yr preprint\endref

\endRefs
\enddocument